\documentclass[14]{article}
\usepackage{latexsym,amssymb,amsmath,amscd,amsthm,amsxtra}
\usepackage{amsmath}
\usepackage{amsfonts}
\usepackage{amssymb}
\usepackage{graphicx}%
\newtheorem{thm}{Theorem}[section]
\newtheorem{prop}[thm]{Proposition} \newtheorem{lem}[thm]{Lemma}
\newtheorem{cor}[thm]{Corollary} \newtheorem{dfn}[thm]{Definition}
 \newtheorem{rmk}[thm]{Remark}
\newtheorem{ex}[thm]{Example} 

\newcommand {\pf}{\noindent{\bf Proof.}\ }
\newcommand{\complex}{{\mathbb C}}

\newcommand{\reals}{{\mathbb R}}
\newcommand{\torus}{{\mathbb T}}
\newcommand{\integers}{{\mathbb Z}}

\newcommand{\cala}{{\cal A}}

\newcommand{\calc}{{\cal C}}
\newcommand{\cald}{{\cal D}}

\newcommand{\calf}{{\cal F}}
\newcommand{\calg}{{\cal G}}

\newcommand{\call}{{\cal L}}

\newcommand{\calq}{{\cal Q}}

\newcommand{\calt}{{\cal T}}

\newcommand{\calu}{{\cal U}}
\newcommand{\calw}{{\cal W}}

\newcommand{\frakk}{\mathfrak{K}}

\title{Deformation Quantization of Pseudo Symplectic(Poisson)
Groupoids}
\author{Xiang Tang \\Department of Mathematics\\ University of California\\ Berkeley, CA
94720 USA\\ {\small(xtang@math.berkeley.edu).}}
\date{}
\begin{document}
\maketitle

\begin{abstract}
We introduce a new kind of groupoid---a pseudo \'etale groupoid,
which provides many interesting examples of noncommutative Poisson
algebras as defined by Block, Getzler, and Xu. Following the idea
that symplectic and Poisson geometries are the semiclassical
limits of the corresponding quantum geometries, we quantize these
noncommutative Poisson manifolds in the framework of deformation
quantization.
\end{abstract}

\begin{center}
{\bf Dedicated to A. Weinstein on his 60th birthday}
\end{center}

\section{Introduction}
About ten years' ago, following an idea from deformation theory,
Block and Getzler in \cite{bl-ge:quantization}, and Xu in
\cite{xu:noncommutative} independently introduced the notion of a
Poisson structure on an associative algebra.

\begin{dfn}
\label{dfn:poisson}A Poisson structure $\Pi$ on an associative
algebra $A$ is an element in the second Hochschild cohomology
$H^2(A; A)$ such that $[\Pi, \Pi]$ is 0, where $[\ ,\ ]$ is the
Gerstenhaber bracket.
\end{dfn}
In \cite{bl-ge:quantization} and \cite{xu:noncommutative}, authors
showed many interesting examples. However, for about ten years,
there were no further results on this subject.

In Section 2, we will introduce a new kind of groupoid, a pseudo
\'etale groupoid, generalizing the notion of an \'etale groupoid
and show that under an extra assumption, the corresponding
groupoid algebras provide a large amount of new examples of
noncommutative Poisson algebras as defined in Definition
\ref{dfn:poisson}.

In physics, a noncommutative Poisson algebra as a phase space of a
noncommutative field theory, and the quantization of this
noncommutative Poisson algebra corresponds to the quantization of
the corresponding field theory. After introducing the
noncommutative Poisson algebras, we study the deformation
quantization of these algebras.

We make the following definition of a formal deformation
quantization of a noncommutative Poisson algebra.
\begin{dfn}
\label{dfn:quant-noncommutative-poisson}Let $(A, \pi)$ be a
noncommutative Poisson algebra,  and let $A[[\hbar]]$ be the
linear space of formal power series with coefficients in $A$. A
{\bf formal deformation quantization} of $(A, \pi)$ is an
associative product $\star \index{$\star$}$ on $A[[\hbar ]]$,
\[
c(\hbar)=a(\hbar)\star \index{$\star$} b(\hbar) \stackrel {def}{=}
\sum _{k=0}^{\infty}\hbar^k c_{k}
\]
satisfying the following properties:
\begin{enumerate}
\item $c_k$ is $\complex[[\hbar]]$ bilinear, \item $c_0=a_0 \cdot
b_0$,
\item
$ a_0\star \index{$\star$} b_0 =c_0 + \hbar i \pi (a_0,
b_0)+o(\hbar). \footnote {$o(\hbar)$ stands for terms with higher
order than $\hbar$. }$
\end{enumerate}
\end{dfn}
\begin{rmk}
The locality assumption of a star product is very important.
Without the locality condition, it is not even known whether the
2nd Hochschild cohomology of the algebra of smooth functions on a
manifold is equal to the space of bivector fields. In Definition
\ref{dfn:quant-noncommutative-poisson}, we have to drop the
locality assumption, since we do not know its corresponding
noncommutative analogue. It would be interesting and useful to
have a notion with which to replace the notion of locality in
noncommutative geometry.
\end{rmk}

In Section 3, we show that the noncommutative Poisson algebra
introduced by a pseudo \'etale groupoid can be formally
deformation quantized. In addition to this, we also discuss
various extension of deformation quantization. Firstly, we
consider existence of traces on the quantized algebras. We show
that if certain modular classes of the corresponding groupoids
vanish, then there is a closed star product as defined by Flato,
Connes, and Sternheimer. Secondly,  instead of taking $\hbar$ as a
formal parameter, we want to look at $\hbar$ as a real number as
in the following definition:
\begin{dfn}
\label{dfn:strict-defor}Let $A$ be a $C^*-$algebra, and
$A_{\infty}$ be its dense $\ast -$subalgebra closed under
holomorphic function calculus, which represents the smooth
algebra. Let $\pi$ be a Poisson structure on $A_{\infty}$. Then
{\bf a strict deformation quantization} of $(A, A_{\infty}, \pi)$
is a set $(\star \index{$\star$} _{\hbar}, ^{*_{\hbar}},
\parallel \cdot \parallel_{\hbar})$ parameterized by $\hbar$, with $\hbar$ in a
closed subset $I$ of the real line containing $0$ as a
non-isolated point. For each $\hbar$, $\star \index{$\star$}
_{\hbar}$ is an associative product, $^{*_{\hbar}}$ is an
involution, and $\parallel \cdot
\parallel_{\hbar}$ is a $C^*-$norm on $A_{\infty}$ satisfying the following
properties:
\begin{enumerate}
\item When $\hbar=0$, the product, involution, and norm come from those on
$A$. \item The completions of $A_{\infty}$ for the various
$C^*-$norms form a continuous field of $C^*-$algebras over $I$.
\item For $a, b\in A_{\infty}$, we have
\[
\parallel (a\star \index{$\star$} _{\hbar} b-a\star \index{$\star$} _0 b -i\hbar\pi (a,b)\parallel_{\hbar} \to
0,\ \ \ \ \hbar \to 0.
\]
\end{enumerate}
\end{dfn}
In 3.4, we discuss strict deformation quantizations in some
interesting examples, which were used in our proof of Morita
invariance of noncommutative tori under the $SO(n,n|\integers)$
action in \cite{tw:dirac}.

The author would like to thank his thesis advisor Professor
Weinstein and coadvisor Professor Rieffel for many helpful
discussions and suggestions.

\section{Noncommutative Poisson algebra} In this section, we
introduce the notion of a pseudo \'etale groupoid, which is a
generalization of an \'etale groupoid. Our main aim is to prove
that under an extra assumption, a pseudo Poisson groupoid defines
a noncommutative Poisson structure on the smooth groupoid
algebra\footnote{See Definition \ref{dfn:smooth-gpd-alg} for a
smooth groupoid algebra. \index{smooth groupoid algebra} }
(Theorem \ref{thm:nc-Poisson-reg} and \ref{thm:coboundary}).
\subsection{Pseudo Poisson groupoid} In this subsection, we introduce
the main object of this paper---a pseudo Poisson groupoid.

\subsubsection{Pseudo \'etale groupoid} A pseudo \'etale groupoid
is a generalization of an \'etale groupoid. The source and target
maps of an \'etale groupoid are defined to be local
diffeomorphisms, so that each element of an \'etale groupoid
defines a local diffeomoprhism on the unit space, from a
neighborhood of an element's source to a neighborhood of its
target. The \'etale assumption of the source and target maps is so
restrictive that a transformation groupoid of a group acting on a
manifold is \'etale if and only if the group is \'etale
(discrete). In the definition of a pseudo \'etale groupoid, we
want to weaken this \'etale assumption, but still keep the
property that each element of the groupoid defines an
infinitesimal diffeomorphism on the unit space, such that the
definition includes all transformation groupoids. A natural way to
achieve this is to choose an infinitesimal bisection at each point
of a groupoid  so that  this choice is compatible with groupoid
operations.

\begin{dfn}
Let $\cald$ be a distribution on a groupoid $\calg$. We say
$\cald$ is {\bf multiplicative} if for any $\alpha,\ \beta,\
\gamma \in \calg$ with $\alpha =\beta \cdot \gamma$ and $X\in
\cald |_{\alpha}$, there are paths $x(t),\ y(t)$ and $z(t)$
satisfying the following properties:
\begin{enumerate}
\item $x(t)=y(t)\cdot z(t)$, and $\dot {x}(t)\in\cald |_{x(t)},\ \dot{y}(t)\in \cald |_{y(t)},\
\dot {z}(t)\in \cald |_{z(t)}$,
\item $x(0)=\alpha,\ y(0)=\beta,\ z(0)=\gamma$, and
$\dot{x}(0)=X$.
\end{enumerate}
\end{dfn}
An example of a multiplicative distribution is the one dimensional
distribution generated by a multiplicative vector field (see
\cite{Mx:multi-vector}) on a Lie groupoid.\footnote{In this sense,
we can look at a multiplicative distribution as a generalization
of a multiplicative vector field.} We will see more interesting
examples later.

\begin{dfn}
An {\bf \'etalification \index{\'etalification}} of a Lie groupoid
$\calg$ is an integrable subbundle $\calf \index{$\calf$}$ of $T
\calg$, satisfying the following conditions:
\begin{enumerate}
    \item \label{compl}$\calf \index{$\calf$}$ is complementary to the $s$ and $t$-fibers. Therefore,
    $s_*$(and $t_*)$ induces an isomorphism between $\calf \index{$\calf$} |_{\gamma}$
and $T_{s(\gamma)}\calg_0$ (and $T_{t(\gamma)}\calg_0$). ($s$ and
$t$ are the source and target maps of $\calg$.)
    \item $\calf \index{$\calf$} |_{\calg^{(0)}}=T \calg ^{(0)}$, where $\calg^{(0)}$ is the unit space of $\calg$.
    \item $\calf \index{$\calf$}$ is a multiplicative distribution.

This assumption together with condition \ref{compl} implies that
    \[
    \dot {y}(0)=(t_{*})^{-1}(X),\ \dot {z}(0)=(s_{*})^{-1}(X).
    \]
\end{enumerate}
\end{dfn}

\begin{dfn}
A {\bf pseudo \'etale groupoid} is a Lie groupoid with a chosen
\'etalification. \index{\'etalification}
\end{dfn}
\begin{rmk}
An \'etalification \index{\'etalification} of a Lie groupoid
redefines the ``topology" on the groupoid, making it behave like
an \'etale groupoid.
\end{rmk}

The following is a list of examples of pseudo \'etale groupoids.
\begin{ex}
\begin{enumerate}
\label{example}
    \item \'Etale groupoid. An \'etale groupoid has  a natural \'etalification\index{\'etalification}, its tangent bundle.
    \item Transformation groupoid $M\rtimes G\rightrightarrows M$. The manifold of
    $M\rtimes G$ is $M\times G$, and the $s,\ t$
maps are defined by $s((x, \gamma))=x,\ t((x, \gamma))=x\cdot
\gamma$ ($\cdot\gamma$\ is the action of $\gamma$ on
    $M$). We choose an \'etalification \index{\'etalification} to be the M-component of the tangent
    bundle $T (M\times G)$.
    \item Pseudo group. In \cite{w1:groupoid}, a pseudo group is defined as ``an algebraic structure
    whose elements consist of selected homeomorphisms between
    open subsets of a space, with the composition of two
    transformations defined on the largest possible domain". The
``germs" of the elements of a pseudo group form a
    groupoid, which we will also call a pseudo group.\footnote{In general, a pseudo group is a not a Lie
    groupoid. Here, we only look at those cases which are Lie groupoids.}
    An \'etalification \index{\'etalification} of a pseudo group $\calg$ can be chosen as follows. A typical element of a pseudo group
    is of the form $(x, \phi, y)$ with $x, y$ two arbitrary elements of
    a manifold $M$ and $\phi$ a germ of a local diffeomorphism from $x$ to $y$.
    For $\gamma=(x, \phi, y) \in \calg$, we choose a representative $\widetilde {\gamma}$ for $\gamma$,
    which defines a diffeomorphism from a neighborhood $U_x$ of $x$ to $U_y$ of $y$, and also induces a
    diffeomorphism from any other point in $U_x$ to its image under
    $\widetilde {\gamma}$. Therefore, $\widetilde {\gamma}$ defines a local bisection near
    $\gamma$. We define an \'etalification \index{\'etalification} of $\calg$ to be the collection of all germs of
    $\widetilde{\gamma}$.
    \item Pair groupoid of $\torus ^k\times \torus ^k \rightrightarrows
    \torus^k$ ($\torus=\reals / \integers$). For simplicity, we describe $\torus \times \torus $, which is represented
    by a unit square $[0,1 ]\times [0,1]$ with edges identified.
    The base of the pair groupoid is the diagonal interval connecting $[0,0]$\ and $[1,1]$. The
    source and target maps are the projections along the $x$ and
    $y$ directions. An \'etalification \index{\'etalification} is chosen to be the subtangent bundle in the diagonal direction.
    \item Germs of bisections on a
groupoid. Generally, a Lie groupoid $\calg$ may not have any
\'etalification \index{\'etalification}. It is not hard to see
that germs of bisections on $\calg$ form a pseudo \'etale groupoid
having a natural \'etalification.\index{\'etalification} An
\'etalification \index{\'etalification} can be chosen in the same
way as a pseudo group in Example 3.
\end{enumerate}
\end{ex}

From the above, we have seen many interesting examples of pseudo
\'etale groupoids. Nevertheless, as was mentioned, not every Lie
groupoid admits an \'etalification.\index{\'etalification} A
counter example is given by the pair groupoid $S^2\times S^2
\rightrightarrows S^2$.

We prove this by contradiction. Suppose that there is an
\'etalification \index{\'etalification} $\calf \index{$\calf$}$ on
the pair groupoid $S^2\times S^2\rightrightarrows S^2$. Choose
$x\in S^2$. At $(x,x)$, we fix one nonzero element $v$ in the
tangent space of the diagonal $S^2$. Since $s_*$ is an isomorphism
between $T_{(x,x)}S^2$ (tangent bundle of the diagonal $S^2$) and
$\calf \index{$\calf$}|_{(y, x)}$, we conclude that
$(s_{*})^{-1}(v)$ induces a nonzero element in the \'etalification
\index{\'etalification} at $(y,x),\ \forall y \in S^2$. Obviously,
at each point $(y, x),\ y\in S^2$, $\ker (s_*)$ and $\ker (t_*)$
are transverse to each other and span the whole tangent space. As
$\calf \index{$\calf$}|_{(y,x)}$ is transverse to both fibers, it
is isomorphic to $\ker (t_*)$ by projection along $\ker (s_*)$. By
this isomorphism, $(s_{*})^{-1}(v)$ is projected onto $\ker
(t_*)$, which is not equal to zero anywhere on the source fiber of
$(y, x )$. This is impossible by the fact that $\ker (t_*)$ is the
tangent bundle of the $t-$fiber at $(x, x)$, which is $T S^2$, and
that $S^2$ has no nowhere vanishing vector fields. Hence, we
conclude that there is no \'etalification \index{\'etalification}
on the pair groupoid $S^2\times S^2\rightrightarrows S^2.$ In
general, by the same type of argument, we can easily prove the
following statement.
\begin{lem}
A pair groupoid $M\times M \rightrightarrows M$ has an
\'etalification \index{\'etalification} if and only if $M$ is
parallelizable.
\end{lem}

The following proposition explains the above lemma.
\begin{prop}
\label{prop:conn-tfiber} Let $(\calg, \calf \index{$\calf$})$ be a
pseudo \'etale groupoid. Since $\calf \index{$\calf$}$ is
transverse to t-fibers (s-fibers), locally $\calf \index{$\calf$}$
induces identifications of t-fibers (s-fibers) along paths in
$\calg^{(0)}$ (an Ehresmann connection), which defines a
connection on $\ker t_*|_{\calg^{(0)}}$ ($\ker
s_*|_{\calg^{(0)}}$). This connection is flat. We call this
connection the canonical connection of $(\calg, \calf
\index{$\calf$})$.
\end{prop}
$\pf$ The existence of this connection is already explained in its
statement. To see that the Ehresmann connection is flat, we notice
that $\calf \index{$\calf$}$ is closed under the Lie bracket and
$t_*$ is a Lie algebra homomorphism of vector fields. $\Box$

There is an another description of a pseudo \'etale groupoid. We
state it without proof\footnote{It is a straightforward check.}.

\begin{prop}
\label{prop:sub-tangent-gpd}An \'etalifiction of a Lie groupoid
$\calg$ is a subtangent groupoid of $T\calg$ which is integrable
and transversal to the source and target fibers of $T\calg$.
\end{prop}

\subsubsection{Poisson structure} Having introduced the notion of a
pseudo \'etale groupoid, we will define and study symplectic
(Poisson) structures on this groupoid.
\begin{dfn}
A symplectic structure on a pseudo \'etale groupoid $(\calg, \calf
\index{$\calf$})$ is a closed 2-form $\omega$ on $\calg$,
satisfying the following conditions:
\begin{enumerate}
\item $\omega | _{\calf \index{$\calf$}}$ makes $\calf \index{$\calf$}$ into a symplectic
bundle; \item $\omega |_{\calf \index{$\calf$}}$ is invariant
under $s^*, t^*$.
\end{enumerate}

A pseudo \'etale groupoid equipped with a symplectic structure is
called a pseudo symplectic groupoid.
\end{dfn}
\begin{dfn}
A Poisson structure on a pseudo \'etale groupoid $(\calg, \calf
\index{$\calf$})$ is a bivector $\pi$ on $\calg$, satisfying the
following conditions:
\begin{enumerate}
\item $\pi |_{\calf \index{$\calf$}}$ makes $\calf \index{$\calf$}$ into a Poisson bundle, which makes $\calf^* \index{$\calf$}$ into a Lie algebroid.
\item $\pi |_{\calf \index{$\calf$}}$ is invariant under $s^*, t^*$.
\end{enumerate}

A pseudo \'etale groupoid equipped with a Poisson structure is
called a pseudo Poisson groupoid.
\end{dfn}
\begin{rmk}
There are two reasons to call these objects ``pseudo". One is that
generally they are not \'etale groupoid, but have similar
properties. The other is that they are not the symplectic
(Poisson) groupoids in the sense of Weinstein.
\end{rmk}

On a symplectic (Poisson) manifold, there is the famous Darboux's
theorem that locally the manifold looks like the standard
symplectic (Poisson\footnote{The Darboux's theorem in Poisson
geometry is more complicated, which should be called
linearization. See \cite{aw:noncom}.}) linear space, which we
recall in the framework of a pseudo \'etale groupoid.
\begin{thm}
Let $(\calg, \calf \index{$\calf$}, \omega)$ be a pseudo
symplectic groupoid. For any $\gamma \in \calg$, there is an
integral submanifold $U$ of $\calf \index{$\calf$}$ through
$\gamma$. In a suitable coordinate system on $U$, $\omega$\ is
expressed by
\[
\sum _{i=1}^{n}dx^i \wedge dy^i .
\]
\end{thm}
\begin{thm}
Let $(\calg, \calf \index{$\calf$}, \pi)$ be a pseudo Poisson
groupoid. For any $\gamma \in \calg$, there is an integral
submanifold $U$ of $\calf \index{$\calf$}$ through $\gamma$. In a
suitable coordinate system of $U$, $\pi$ is expressed by
\[
\sum _{i}\frac{\partial}{\partial q_i}\wedge \frac
{\partial}{\partial p_i}+\frac {1}{2}\sum _{i,j}\phi _{ij}\frac
{\partial}{\partial y_i}\wedge \frac {\partial}{\partial y_j}, \ \
\ \ \ \ \phi _{ij}(0)=0.
\]
\end{thm}

\begin{ex}{\bf of pseudo symplectic (Poisson) groupoids \index{pseudo symplectic (Poisson) groupoid}}
\begin{enumerate}
\label{ex:poisson} \item Transformation groupoid. A transformation
groupoid with an invariant symplectic (Poisson) structure on the
manifold. \item \'Etale groupoid. An \'etale groupoid with an
invariant symplectic (Poisson) structure. \item Pseudo group. For
a symplectic (Poisson) manifold, if we ask the local
diffeomorphism to be symplectic (Poisson), the corresponding
pseudo group forms a pseudo symplectic (Poisson) groupoid
\index{pseudo symplectic (Poisson) groupoid}. \item Pair groupoid
$\torus^k \times \torus ^k\rightrightarrows \torus^k$. For a pair
groupoid $\torus^k \times \torus^k \rightrightarrows \torus^k$, we
associate a constant symplectic (Poisson) structure on $\torus^k$,
and $s^*$ induces a symplectic (Poisson) structure on the entire
groupoid. In this way, the pair groupoid forms a pseudo symplectic
(Poisson) groupoid \index{pseudo symplectic (Poisson) groupoid}.
\item Orientable contact manifold $(M, \alpha)$. On a contact
manifold, the Reeb vector field $R$ generates a 1-dimensional
foliation on $M$. We choose a submanifold of $M$ which is
transverse to the foliation. It is easy to check that $d\alpha$
defines a symplectic form on the transversal, which is invariant
along the foliation, and the reduced foliation groupoid to this
transversal forms a pseudo symplectic groupoid. \item
\label{ex:dirac} Dirac manifold $(M, L)$ with a constant rank
characteristic distribution. The characteristic distribution
$L\cap T M$ of a Dirac manifold, which is integrable, forms a
foliation on M. As in the contact case, we choose a submanifold of
$M$ which is transverse to the foliation. On this transversal,
there is a natural Poisson structure invariant along the
foliation. This makes the reduced foliation groupoid on the chosen
transversal into a pseudo Poisson groupoid.
\end{enumerate}
\end{ex}
The last two examples are specific examples of \'etale groupoids
with an invariant symplectic (Poisson) structure. We point them
out to show a wide applicability of our results.

Given a pseudo symplectic groupoid $\calg$, its \'etalification
\index{\'etalification} $\calf \index{$\calf$}$ forms a symplectic
vector bundle on $\calg$. So we can ask for a symplectic
connection on the bundle.

\begin{dfn}
\label{dfn:connection} A pseudo symplectic connection on a pseudo
symplectic groupoid is a symplectic connection on the symplectic
bundle $(\calf \index{$\calf$} | _{I}, \omega)$ invariant under
$s,t$, where $I$ denotes the leaves of $\calf \index{$\calf$}$.
\end{dfn}
\begin{lem}
If a pseudo symplectic groupoid $\calg$ is proper, then there is a
pseudo symplectic connection on it.
\end{lem}
$\pf$ A symplectic connection on $\calf \index{$\calf$}$ always
exists by the same arguments as in Proposition 2.5.2 in
\cite{Fe:deformation}. So we only have to calculate its invariance
under $s, t$. When $\calg$ is proper, one can choose a Haar system
\index{Haar system} (see 2.1.3) on $\calg$, and integrate the
connection by this Haar system. It is easy to check that the
integrated connection is a pseudo symplectic connection. $\Box$

\begin{rmk}
A similar result holds for the regular Poisson case. On a general
Poisson manifold, there is usually no Poisson connection. It can
be easily shown that a Poisson manifold has a Poisson connection
if and only if the Poisson structure is regular.
\end{rmk}

\subsubsection{Groupoid algebra} The smooth groupoid algebra
\index{smooth groupoid algebra} of a Lie groupoid consists of
smooth functions on the groupoid with a convolution product, which
is an important example of noncommutative differentiable
manifolds. In the rest of this section, we plan to translate a
symplectic (Poisson) structure on a pseudo \'etale groupoid to a
noncommutative Poisson structure on the corresponding smooth
groupoid algebra\index{smooth groupoid algebra}. In this
subsubsection, we will recall the concept of a smooth groupoid
algebra\index{smooth groupoid algebra}. There are two issues we
will talk about. One is how to define a convolution product, the
other is the definition of a smooth function on a Lie groupoid.

There are usually two ways to define a convolution product on a
groupoid algebra. One is to use a Haar system,\index{Haar system}
the other is to replace functions by half densities. The two
methods yield the same algebra. For the details and the
equivalence of the two versions, readers are referred to
Paterson's book \cite{P1:groupoid}. In this thesis, we will define
a groupoid algebra by a Haar system\index{Haar system}. We begin
by recalling the definition of a Haar system\index{Haar system}.

(\cite{P1:groupoid}, Definition 2.3.2) A smooth left Haar system
\index{Haar system} for a Lie groupoid $\calg$ is a family $\{
\lambda ^u\}, u\in \calg ^{(0)}$ (the unit space of $\calg$),
where each $\lambda ^u$ is a positive regular Borel measure on the
manifold $\calg ^u$ (the t-fiber of $\calg$ at unit $u$.), such
that the following axioms are satisfied:
\begin{enumerate}
    \item If $(V, \psi)$ is a t-fiberwise product open subset of
    $\calg,\ i.e.\ V\cong t(V)\times W$, and if  $\lambda _W$ is
    Lebesgue measure on $\reals ^k$ restricted to $W$, then for
    each $u\in t(V)$, the measure $\lambda ^u \cdot \psi _u$ is
    equivalent to $\lambda _W$, and the map $(u,\ w)\to d(\lambda ^u \cdot \psi _u)/ d \lambda _W
    (w)$ belongs to $C^{\infty}(t(V)\times W)$ and is strictly
    positive;
    \item For any $x\in \calg$ and $f\in C_c (\calg)$, we have
    \[ \int _{\calg ^{s(x)}}f(xz)d\lambda ^{s(x)}(z)= \int _{\calg
    ^{t(x)}} f(y)d\lambda ^{t(x)}(y). \]
\end{enumerate}
From Theorem 2.3.1 of \cite{P1:groupoid}, we know the smooth left
Haar system \index{Haar system} on a Lie groupoid always exists.

In the case of a pseudo \'etale groupoid, we have a special
requirement on the choice of a Haar system\index{Haar system}. As
is suggested in \cite{Mx:multi-vector}, a multiplicative vector
field defines a derivation on a smooth groupoid
algebra\index{smooth groupoid algebra}. Accordingly, we want every
section of an \'etalification\index{\'etalification} to define a
derivation (see Lemma \ref{derivation}). However, this is only
true when using the half density definition of a groupoid algebra.
When using a Haar system\index{Haar system}, we usually get three
terms out of differentiating a product of two functions. Two of
them are the expected terms from the Leibniz rule, but the other
contains a derivative respect to the chosen Haar system
\index{Haar system}. To get rid of this term in the
differentiation, we make the following definition.

\begin{dfn}
\label{dfn:trans-haar} Let $(\calg, \calf \index{$\calf$})$ be a
pseudo \'etale groupoid.  According to the definition, each
t-fiber of $\calg$ is transversal to the leaves of $\calf
\index{$\calf$}$. We call a Haar system transversal if it defines
a transversal measure to the foliation of $\calf \index{$\calf$}$.
We call this kind of Haar system \index{Haar system} a {\bf
transversal Haar system}\index{transversal Haar system}.
\end{dfn}

For a pseudo \'etale groupoid, a transversal Haar system
\index{Haar system} does not always exist.  Its existence is
determined by the following property of the canonical connection
defined in Proposition \ref{prop:conn-tfiber}.

\begin{prop}
\label{prop:exit-trans-measure}For any pseudo \'etale groupoid
$(\calg,\ \calf \index{$\calf$})$, if the canonical connection
defined in Proposition \ref{prop:conn-tfiber} has holonomy in
$SL(k; \reals) $ $(k=dim(\ker t_*))$,  there is a transversal Haar
system \index{Haar system}.
\end{prop}
$\pf$ Consider the density bundle $\hat {H}$ of $ker t_* | _{\calg
^{(0)}}$ on $\calg^{(0)}$, which is flat. It is easy to see that
there is a one to one correspondence between the Haar systems
\index{Haar system} and the  sections of $\hat {H}$. Given a
section $\lambda \in \Gamma (\hat {H})$, we can check that
$\Lambda(\alpha)\stackrel{def}{=}(\alpha
^{-1})^{*}(\lambda(s(\alpha)))$ ($\alpha: s(\alpha)\mapsto
\alpha=\alpha\cdot s(\alpha)$ while $\alpha ^{-1}:\alpha \mapsto
s(\alpha)=\alpha ^{-1}\cdot \alpha$) defines a left invariant
volume form along the $t-$fibers, which is a Haar
system\index{Haar system}. For the other direction, a smooth Haar
system \index{Haar system} when restricted to the unit space
defines a section of $\hat {H}$. Obviously, the two maps are
inverse to each other.

Since the holonomy of the canonical connection is in $SL(k;
\reals)$, $\hat {H}$ is trivial and has global flat sections. We
choose one global flat section,  which is everywhere nonvanishing,
denoted by $\lambda$. We show that $\Lambda(\alpha)
\stackrel{def}{=}(\alpha ^*)^{-1}(\lambda(s(\alpha)))$ is a
transversal Haar system \index{Haar system}.

For $X\in \calf \index{$\calf$}$, locally we have
\[
\begin{array}{ll}
\call_{X}\Lambda&=s^*\circ(s^{-1})^*(\call_{X}\Lambda)\\
&=s^*(\call_{s_*(X)}(s^{-1})^*(\Lambda))\\
&=s^*(\call_{s_*(X)}(\lambda))\\
&=0
\end{array}
\]
where $\alpha ^{-1}:\alpha \mapsto s(\alpha)$ is equal to $s$.
$\Box$

There are many examples satisfying the condition of Proposition
\ref{prop:exit-trans-measure}:
\begin{ex}
\begin{enumerate}
\item \'Etale groupiod. In the case of an \'etale groupoid, $k=dim(\ker
t_*)=0$. The bundle $\ker(t_*)|_{\calg^{(0)}}$ is always trivial.
\item Transformation groupoid. $M\times G$ is a trivial $G$-bundle
over $M$, which makes $\ker (t_*)|_{M}$ also a trivial bundle.
\item Pseudo groupoid. If we require homomorphisms in a
pseudo group to preserve a metric, then the honolomy of the
canonical connection is contained in $SO(k; \reals)$ which is in
$SL(k; \mathbb {R})$.
\item Pair groupoid $\torus ^n\times \torus^n\rightrightarrows
\torus^n$. It is straight forward to check that the holonomy of
the canonical connection is also trivial.
\end{enumerate}
\end{ex}

In the following, we will always assume that there exists a
transversal Haar system \index{Haar system} on a pseudo \'etale
groupoid, and define a convolution product on $C_c
^{\infty}(\calg)$,
\[
f\diamond \index{$\diamond$} g(\alpha)\stackrel {def}{=}\int
_{\alpha=\beta \cdot \gamma}f(\beta)g(\gamma)d \lambda.
\]
It is not hard to check that $\diamond \index{$\diamond$}$ is
associative.

Now we discuss the definition of a smooth function on a Lie
groupoid. The problem comes from the possibility that a Lie
groupoid may not be a Hausdorff manifold. When a Lie groupoid is
Hausdorff, compactly supported smooth functions on it form an
associative algebra under the convolution product. But in a
non-Hausdorff case, the sum of two compactly supported smooth
functions may not be smooth. Therefore, we have to enlarge the set
of compactly supported smooth functions to make it closed under
summation and convolution. There are two definitions of a smooth
groupoid algebra \index{smooth groupoid algebra}. One can be found
in Connes' book \cite{C:nc}, the other is used in Crainic and
Moerkijk's work (see \cite{crainic:cyclic}, \cite{cm:homology}).
The two definitions agree when the groupoid is Hausdorff, in
particular when it is proper. In \cite{C:nc}, for a non-Hausdorff
smooth groupoid $\calg$, $C_c ^{\infty}(\calg)$ is defined to be
the complex linear subspace of all functions on $\calg^{(1)}$
generated by $C_c ^{\infty}(U)$ for all coordinate domains $U$.
The smooth groupoid algebra \index{smooth groupoid algebra}
$A(\calg)$ is the smallest subalgebra of the algebra of compactly
supported Borel functions on $\calg$ which is closed under the
convolution product $\diamond \index{$\diamond$}$ and contains
$C_c ^{\infty}(\calg)$. However, even though this definition is
the ``smallest" one, there is a ``better" algebra introduced in
\cite{crainic:cyclic} and \cite{cm:homology}, whose cyclic and
Hochschild homology are easier to calculate. In the rest of this
subsection, we will recall Crainic and Moerdijk's definition of a
smooth groupoid algebra \index{smooth groupoid algebra}.

\begin{dfn}
\label{dfn:gamma-c}(2.11, \cite{crainic:cyclic}) Let $X$ be a
topological space. A sheaf $\cala \in Sh(X)$ is called $c-$soft if
there is a Hausdorff open covering $\calu$ of $X$ such that $\cala
|_U \in Sh(U)$ is $c-$soft (see \cite{crainic:cyclic} and
\cite{cm:homology}) for all $U\in \calu$. In this case, define
$\Gamma _c (U; \cala)$ as the image of the map
\[
\bigoplus _U \Gamma _c(U; \cala)\to \Gamma (X_{dis}; \cala),
\]
where $\Gamma(X_{dis};\cala)=\{u:X\to \bigsqcup _{x\in X}\cala_x:
u(x)\in \cala_x,\ \forall x\in \cala_x \}$ ($X_{dis}$ is X,
considered with the discrete topology) and $\Gamma _c(U; \cala)\to
\Gamma (X_{dis}; \cala),\ s\mapsto \bar{s}$ is given by
\[
\bar{s}(x)=germ_x(s)\ for\ x\in U,\ and\ 0\ otherwise.
\]
\end{dfn}
\begin{dfn}
\label{dfn:smooth-fun}(2.14, \cite{crainic:cyclic}) If $M$ is a
manifold, not necessarily Hausdorff, we define $C _c
^{\infty}(M)\stackrel{def}{=}\Gamma _c (M; \calc _M ^{\infty})$
where $\calc _M^{\infty}$ is the sheaf of smooth functions on $M$.
From the Mayer-Vietoris sequence, we have an alternative
description of $C_c ^{\infty}(M)$, as the cokernel of:
\[
\bigoplus _{U, V} C _c ^{\infty}(U\cap V)\to \bigoplus _U C_c
^{\infty}(U)\ \ \ \ (U\in \calu),
\]
where $\calu$ is a Hausdorff open covering of $M$.
\end{dfn}

\begin{dfn}
\label{dfn:smooth-gpd-alg}We define the smooth groupoid algebra
\index{smooth groupoid algebra} of $\calg$ to be $C _c
^{\infty}(\calg)$ defined by Definition \ref{dfn:smooth-fun}.
\end{dfn}

If $\calg$ is pseudo \'etale, the \'etalification
\index{\'etalification} makes groupoid operations on $\calg$ into
local diffeomorphisms, i.e. the left multiplication of $\alpha \in
\calg$ on $\beta \in \calg$ naturally identifies $\beta$'s s-fiber
with $(\alpha \cdot \beta)'s$ s-fiber, and simultaneously the
\'etalification \index{\'etalification} maps $\calf
\index{$\calf$}|_{\beta}$ to $\calf
\index{$\calf$}|_{\alpha\cdot\beta}$ isomorphically. Hence, the
left multiplication of $\alpha$ defines a local diffeomorphism.
With this observation, the same argument as in \cite{bn:cyclic}
shows that $C _c ^{\infty}(\calg)$ is closed under the convolution
product $\diamond \index{$\diamond$}$. This algebra is what we
will work with in the next subsection.

\subsection{Noncommutative Poisson Structure} For a pseudo symplectic (Poisson) groupoid
\index{pseudo symplectic (Poisson) groupoid}, we define
\begin{equation}
\label{Poisson} \Pi (f, g)(\alpha )\stackrel {def}{=}\int _{\beta
\gamma =\alpha}\pi (\alpha)((t^*)^{-1}(df(\beta)), (s^*)^{-1}(dg
(\gamma))),\ \ \ \ \ \forall f, g\in C_c ^{\infty}(\calg),
\end{equation}
where $(t^*)^{-1}(df(\beta))$\ \ (or $(s^*)^{-1}(dg (\gamma))$ )
means that we first restrict $f$\ $(or\ g)$ to the integrated
submanifold near $\beta\ (or\ \gamma)$, then calculate the
differential there to form an element in $\calf \index{$\calf$}
^{*}|_{\beta}\ (or\ \calf \index{$\calf$} ^*|_{\gamma})$, and
finally push the element forward to $\calf \index{$\calf$}
^*|_{\alpha}$ by $(t^*)^{-1}\ (or\ (s^*)^{-1})$. The main
objective of this section is to see when formula (\ref{Poisson})
defines a noncommutative Poisson structure on the groupoid algebra
$C_c^{\infty} (\calg)$.
\subsubsection{Regular Poisson case} The
main result of this subsection is the following theorem, which
generalizes Proposition 2.3 in \cite{bl-ge:quantization}.
\begin{thm}
\label{thm:nc-Poisson-reg} A pseudo regular Poisson\footnote{Here,
``regular" means that $\pi$ has constant rank.} (e.g. pseudo
symplectic ) groupoid $(\calg, \calf \index{$\calf$}, \omega)$ (or
$(\calg, \calf \index{$\calf$}, \pi)$) with a given invariant
connection $\nabla$ naturally defines a noncommutative Poisson
structure on $C_c^{\infty}(\calg)$.
\end{thm}
\begin{rmk}
Here, We will use the invariant connection $\nabla$ to construct a
``coboundary" of  $[\Pi, \Pi]$. We believe that the existence of a
Poisson structure on a groupoid algebra should imply the existence
of a kind of invariant ``connection" in a generalized sense.
However, we do not know how to define this.
\end{rmk}
$\pf$ To prove that $\Pi$ is a noncommutative Poisson structure,
we have to show that $\Pi$ defines a Hochschild 2-cocycle with
$[\Pi,\ \Pi]$ being a 3-coboundary. We start with proving the
following lemma.

\begin{lem}
\label{derivation} When restricted to the leaves of $\calf
\index{$\calf$}$, $d$ has the following formula: $\forall f,\ g,\
\in C_c ^{\infty}(\calg)$, :
\[
d(f\diamond \index{$\diamond$} g)(\alpha )=\int_{\alpha=\beta
\gamma} (t^*)^{-1}(df(\beta))g(\gamma) + \int_{\alpha=\beta
\gamma}f(\beta) (s^*)^{-1}(dg (\gamma)),
\]
where $\diamond \index{$\diamond$}$ stands for the convolution of
groupoid algebra.
\end{lem}
$\pf$ We have to use the multiplicativity of $\calf
\index{$\calf$}$. Suppose that $x\in \calf \index{$\calf$}
|_{\alpha}$. Then the multiplicative assumption on $\calf
\index{$\calf$}$ provides $x(t), y(t), z(t)$, satisfying
\begin{enumerate} \item $x(0)=\alpha,\ y(0)=\beta,\ z(0)=\gamma$,
\item $\dot {x}(t)\in \calf \index{$\calf$} |_{x(t)},\ \dot{y}(t)\in \calf \index{$\calf$} |_{y(t)},\
\dot {z}(t)\in \calf \index{$\calf$} |_{z(t)}$, \item
$\dot{x}(0)=x,\ \dot {y}(0)=(t_*)^{-1}(x),\ \dot
{z}(0)=(s_*)^{-1}(x)$.
\end{enumerate}

 Therefore, $x(f\diamond \index{$\diamond$} g)(\alpha)=$
 \[
\begin{array}{ll}
=&\frac {d}{dt}|_{t=0}f\diamond \index{$\diamond$} g(x(t))\\
=&\frac {d}{dt}|_{t=0}\int _{y(t)z(t) =x(t)}f(y(t))g(z(t))\\
=&\int _{y(t)z(t) =x(t)} \frac {d}{dt}|_{t=0}(f(y(t))g(z(t)))\\
=&\int _{y(t)z(t) =x(t)} (\frac
{d}{dt}|_{t=0}f(y(t)))g(\gamma)+f(\beta)(\frac
{d}{dt}|_{t=0}g(z(t)))\\
=&\int _{\alpha=\beta \gamma}
((t_{*})^{-1}x)f(\beta)g(\gamma)+f(\beta)((s_{*})^{-1}x)g(\gamma)\\
=&\int _{\alpha=\beta
\gamma}(t^*)^{-1}(df(\beta))(x)g(\gamma)+f(\beta)(s^{*})^{-1}(dg
(\gamma)),
\end{array}
\]
where in the third equality we have used the fact that our Haar
system \index{Haar system} is transversal, otherwise there would
be one more term like $f(\gamma)g(\beta)\frac{d}{dt}(\lambda
^{r(x(t))})$.
\begin{lem}
\label{2cocycle} $\Pi$ satisfies the cycle condition,
\[
f\diamond \index{$\diamond$} \Pi (g, h)-\Pi (f\diamond
\index{$\diamond$} g, h)+ \Pi (f, g\diamond \index{$\diamond$}
h)-\Pi (f, g)\diamond \index{$\diamond$} h=0,\ \ \ \ \ \ \forall
f, g, h \in C_c ^{\infty}(\calg).
\]
\end{lem}
$\pf$ By Lemma \ref{derivation}, have
\[
\begin{array}{ll}
\Pi (f\diamond \index{$\diamond$} g, h)(\alpha)&= \int _{\beta
\gamma=\alpha} \pi (\alpha)((t ^*)^{-1}(d(f\diamond
\index{$\diamond$} g)(\beta)), (s^*)^{-1}(dh
(\gamma)))\\
&=\int _{\beta \gamma =\alpha}\pi (\alpha)((t
^*)^{-1}(\int_{\xi \eta=\beta} (t^*)^{-1}(df(\xi))g(\eta)\\
&+ \int_{\xi \eta=\beta}f(\xi) (s^*)^{-1}(dg (\eta))),
(s^*)^{-1}(dh(\gamma))\\
&=\int _{\beta \gamma =\alpha}\int_{\xi \eta=\beta}\pi (\alpha)((t
^*)^{-1}((t^*)^{-1}(df(\xi))g(\eta)\\
&+f(\xi) (s^*)^{-1}(dg (\eta))), (s^*)^{-1}(dh (\gamma)))\\
&=f\diamond \index{$\diamond$} (\Pi(g, h))(\alpha)+\int _{\beta
\gamma =\alpha}\int_{\xi \eta=\beta}\pi
(\alpha)((t^*)^{-1}(df(\xi))g(\eta), (s^*)^{-1}(dh (\gamma))).
\end{array}
\]
Similarly,
\[
\begin{array}{ll}
\Pi (f,g\diamond \index{$\diamond$} h)(\alpha)&=\Pi(f, g)\diamond
\index{$\diamond$} h(\alpha)+\int _{\beta \gamma =\alpha}\int_{\xi
\eta=\beta}\pi (\alpha)((t^*)^{-1}(df(\xi))g(\eta), (s^*)^{-1}(dh
(\gamma))).\ \ \ \Box
\end{array}
\]

To prove Theorem \ref{thm:nc-Poisson-reg}, we still need to show
that $[\Pi, \Pi]$ is a 3-coboundary. To prove this, we use an
invariant Poisson connection $\nabla$ on $(\calg,\calf
\index{$\calf$})$ to define $P _2$ as follows:
\[
P _2(f, g)(\alpha)\stackrel {def}{=}\int _{\alpha =\beta
\gamma}<\pi \otimes \pi (\alpha), (t^*)^{-1}(\nabla
^2f(\beta))\otimes (s^*)^{-1}(\nabla ^2g(\gamma))
>,
\]
where we restrict $f, g$ to leaves of $\calf \index{$\calf$}$ to
construct $\nabla ^2f, \nabla ^2 g$, and pair the tensor $\pi
\otimes \pi$ with an element $\alpha \otimes \beta \in S^2 (T^* M)
\otimes S^2 (T^* M)$ by the following formula
\[
\pi ^{ij}\pi ^{kl}\alpha _{ik}\beta _{jl}.
\]
It is easy to see that $P_2$ is skew-symmetric.
\begin{lem}
\[
\delta P_2 + [\Pi, \Pi]=0.
\]
\end{lem}
$\pf$ Since $\nabla$ is an invariant Poisson connection, we have
\[
\nabla \pi=0,
\]
and
\[
d<\pi ,(t^*)^{-1}df\wedge (s^*)^{-1}dg>=<\pi, (t^*)^{-1}(\nabla ^2
f(\alpha))\otimes (s^*)^{-1}dg+ (t^*)^{-1}df\otimes
(s^*)^{-1}(\nabla ^2 g(\beta))> .
\]
Using the above formulas, we calculate
\[
\begin{array}{ll}
\Pi(\Pi (f, g), h)(\alpha)&= \int _{\alpha=\beta \gamma}\pi
(\alpha)((t^*)^{-1}(d(\Pi (f, g)))(\beta),
(s^*)^{-1}(h)(\gamma))\\
&=\int _{\alpha=\beta \gamma}\int _{\beta=\xi \eta}\pi (\alpha)(
(t^*)^{-1}(d<\pi (\alpha),(t^*)(df)(\xi)\wedge
(s^*)^{-1}(dg(\eta))>),\\
&(s^*)^{-1}(h)(\gamma))\\
&=\int _{\alpha=\beta \gamma}\int _{\beta=\xi \eta}<\pi \otimes
\pi, (t^{*})^{-1}(\nabla ^2f (\xi))\otimes (s^*)^{-1}(dg
(\eta))\otimes (s^*)^{-1}(dh(\gamma))\\
&+(t^*)^{-1}(df(\xi))\otimes (s^{*})^{-1}(\nabla ^2 g
(\eta))\otimes(s^*)^{-1}(dh(\gamma))>;\\
\Pi (f, \Pi (g, h))(\alpha)&=\int _{\alpha=\xi \beta}\int
_{\beta=\eta\gamma \ }<\pi \otimes \pi, (t^{*})^{-1}(df
(\xi))\otimes (s^*)^{-1}(\nabla ^2g
(\eta))\otimes (s^*)^{-1}(dh(\gamma))\\
&+(t^*)^{-1}(df(\xi))\otimes (s^{*})^{-1}(d g
(\eta))\otimes(s^*)^{-1}(\nabla ^2h(\gamma))>.
\end{array}
\]
On the other hand, using the formula $\nabla ^2 (f \diamond
\index{$\diamond$} g)(\alpha)=$
\[
\int_{\alpha =\beta \gamma}(t^{*})^{-1}(\nabla ^2f
(\beta))g(\gamma)+2 (t^*)^{-1}(df
(\beta))(s^*)^{-1}(dg(\gamma))+f(\beta)(s^*)^{-1}(\nabla ^2
g(\gamma)),
\]
we can show
\[
\begin{array}{ll}
&f\diamond \index{$\diamond$} P_2 (g, h)(\alpha)-P_2(f\diamond \index{$\diamond$} g, h)(\alpha)\\
&=\int _{\alpha=\beta\gamma}\int _{\gamma =\xi \eta}<\pi \otimes
\pi (\alpha),
(t^{*})^{-1}(\nabla ^2 f(\beta))g(\xi)(s^*)^{-1}(dh(\eta))\\
&+(t^*)^{-1}(df
(\beta))(t^*)^{-1}(dg(\xi))(s^*)^{-1}(dh (\eta))>,\\
\end{array}
\]
and also,
\[
\begin{array}{ll}
&-P_2(f, g\diamond \index{$\diamond$} h)(\alpha)+P_2(f, g)\diamond \index{$\diamond$} h (\alpha)\\
&=\int _{\alpha =\beta\gamma}\int _{\gamma =\xi \eta} <\pi \otimes
\pi (\alpha), (t^*)^{-1}(df (\beta))\otimes (t^*)^{-1}(dg
(\xi))\otimes (s^*)^{-1} (dh (\eta))\\
&+(t^*)^{-1}(df (\beta))g(\xi)\otimes (s^*)^{-1}(\nabla
^2(dh(\eta)))>.
\end{array}
\]
The equality we need for this lemma easily follows from the above
calculation. $\Box$

We have finished the proof of Theorem \ref{thm:nc-Poisson-reg}. We
see that the above proof strongly depends on the existence of an
invariant connection. Hence it cannot deal with general pseudo
Poisson groupoids. We will look at the Poisson cases in the next
subsection.

\subsubsection{General Poisson case} In \cite{Ko:deformation},
Kontsevich proved his famous formality \index{formality} theorem.
As a corollary, he showed that every Poisson manifold can be
deformation quantized. In this subsection, we discuss an
``equivariant formality \index{formality} theorem" with a groupoid
action.

Let $(\calg, \calf \index{$\calf$})$ be a proper pseudo \'etale
groupoid. On the unit space $\calg^{(0)}$ of $\calg$, a
differentiable manifold, there are two differential graded
algebras which are defined as follows:

\begin{enumerate}
\item Let $\calt^* _{poly}$ be the space of multi-vector
fields on $\calg ^{(0)}$. $(\calt^* _{poly},\ 0,\ [\ ,\ ]\ )$
defines a differential graded Lie algebra (DGLA),  where $0$ is
the $0$ differential, and $[\ ,\ ]$ is the Schouten-Nijenhuis
bracket.
\item Let $\cald^*_{poly}  \index{$\cald^*_{poly}$}$ be the subspace of $Hom_{\complex}
((C^{\infty}(\calg ^{(0)}))^{\otimes *}, C^{\infty}(\calg^{(0)}))$
consisting of multi-differential operators. $(\cald^*_{poly}
\index{$\cald^*_{poly}$},\ \delta ,\ [\ ,\ ]\ )$ also defines a
DGLA, where $\delta$ is the Hochschild differential and $[\ ,\ ]$
is the Gerstenhaber bracket.
\end{enumerate}

Kontsevich, in \cite{Ko:deformation}, proved a formality
\index{formality} theorem that the above two DGLAs are
quasi-isomorphic. When a pseudo \'etale groupoid $\calg$ is
proper, we can consider the corresponding $\calg$ invariant sub
DGLA. By integration, it is easy to check that Kontsevich's proof
of formality \index{formality} theorem still works in the
invariant case. So we have the following theorem.
\begin{thm}
\label{eq-formality}If\ $\calg$ is proper, $((\calt^*_{poly}
\index{$\calt^*_{poly}$})^{\calg},\ 0,\ [\ ,\ ]\ ) $ is
quasi-isomorphic to $((\cald^*_{poly}
\index{$\cald^*_{poly}$})^{\calg},\ \delta ,$ $[\ ,\ ])$, where
$(\calt ^*_{poly})^{\calg}$ and $(\cald^*_{poly}
\index{$\cald^*_{poly}$})^{\calg}$ are the corresponding $\calg$
invariant subspace.
\end{thm}
Proof: Following the proof of Theorem 4.6.2 in
\cite{Ko:deformation}, we globalize the local formality
\index{formality} theorem on $\reals ^n$ by the following steps
\[
\begin{array}{ll}\calt _{poly}(\calg^{(0)})[1]_{formal}&\stackrel {1}{\to}\Gamma
(\calt_{s^{aff}}\to {T[1]\calg ^{(0)})}_{formal}\\
&\stackrel{2}{\to} \Gamma (\cald _{s^{aff}}\to T[1]\calg
^{(0)})_{formal}\stackrel{3}{\leftarrow} \cald
_{poly}{(\calg^{(0)})[1]}_{formal}.
\end{array}
\]

In the above maps, 2 is a fiberwise map which is invariant under
$\calg$ action. We have to integrate maps 1 and 3 by $\calg$ to
make them $\calg$ invariant, where we use the properness
assumption of $\calg$. In this way, we obtain a quasi-isomorphism
between the invariant parts of the two $L_{\infty}$ algebras.
$\Box$

\begin{thm}
\label{thm:coboundary} Let $(\calg, \calf \index{$\calf$}, \pi)$
be a proper pseudo Poisson groupoid. $\Pi$ in (\ref{Poisson})
defines a Poisson structure on $(C_c^{\infty}(\calg),\ \diamond
\index{$\diamond$} )$.
\end{thm}
$\pf$ The different part of the proof from Theorem
\ref{thm:nc-Poisson-reg} is to show that $[\Pi, \Pi]$ is a
3-coboundary. We will use the above equivariant formality
\index{formality} Theorem \ref{eq-formality} to construct a
2-cochain whose coboundary is $[\Pi, \Pi]$.

We look at $\pi$ restricted to the unit space $\calg^{(0)}$. As
$\pi$ is a Poisson structure,  there is a Hochschild 2-cocycle
$P_2$ on $\calg^{(0)}$ with
\begin{equation}
\label{coboundary} [\pi, \pi] =\delta P_2.
\end{equation}
Furthermore, because $\pi$ is invariant under $\calg$ and $\calg$
is proper, we can integrate $P_2$ along $\calg$ orbits to make it
invariant.

Then because $s,\ t$ defines a local diffeomorphism between the
leaves of $\calf \index{$\calf$}$ and $\calg _0$, we can pull back
$P_2$ to the entire $\calg$, which is still denoted by $P_2$. We
define
\[
\hat {P}_2(f, g)(\alpha)\stackrel {def}{=}\int _{\alpha
=\beta\gamma}P_2(\alpha)((t^*)^{-1}f(\beta), (s^*)^{-1}g(\gamma)),
\]
where $f,g$ are smooth compactly supported functions on $\calg$.

From equality (\ref{coboundary}), it is straightforward to check
that
\[
\delta \hat {P}_2=[\Pi, \Pi] .\ \ \ \ \ \ \Box
\]

\begin{rmk}
\label{rmk:equ-formality}The assumption of Theorem
\ref{eq-formality} can be weakened. Dolgushev, in
\cite{Do:formality} showed a new way to prove the global formality
\index{formality} theorem from Kontsevich's local formula by a
Fedosov (see \cite{Fe:deformation}) type resolution. Dolgushev
proved that, with a connection, there is a global
quasi-isomorphism between the two $L_{\infty}-$algebras. To get an
equivariant or invariant quasi-isomorphism, we only need to
require that there is an invariant connection. Following his
methods without repeating the details, we can easily show that if
there is an invariant connection, then $\Pi$ defines a Poisson
structure on $C_c^{\infty}(\calg)$.
\end{rmk}
\begin{rmk}
If $\Pi$ is the first term in a deformation quantization of a
noncommutative algebra, then we can show that $[\Pi, \Pi]$ must be
a 3-coboundary as a Hochschild cocycle. This statement is closely
related to Theorem \ref{eq-formality}.
\end{rmk}

In conclusion, in this section we have introduced a new class of
Lie groupoid---pseudo symplectic (Poisson) groupoid \index{pseudo
symplectic (Poisson) groupoid}, which provides many new examples
of noncommutative Poisson algebras. This is the object we will
quantize in the next chapter.

\section{Quantization of a Pseudo \'Etale Groupoid} In this
section, we consider quantization of pseudo \'etale groupoids. In
3.1, we construct a star product on a noncommutative Poisson
manifold introduced in Theorem \ref{thm:nc-Poisson-reg} and
\ref{thm:coboundary}. In 3.2 and 3.4, we discuss two extensions of
quantization problems, closed star products and strict deformation
quantizations. In 3.3, we demonstrate our quantization methods on
a typical example, a transformation groupoid.

\subsection{Quantization} In this section, using
Fedosov's machinery, we construct a formal deformation
quantization of the noncommutative Poisson algebra defined in the
last section.

If we forgot about the groupoid structure, a pseudo Poisson
groupoid is a Poisson manifold. We know how to quantize a Poisson
manifold from \cite{Ko:deformation} and \cite{Cf1:local-global}.
However, the quantization of this manifold is of less interest,
since it does not contain the information of the groupoid
operations. Instead, we should consider a quantization of the
noncommutative Poisson algebra, which is more complicated than the
commutative Poisson algebra of smooth functions on a Poisson
manifold. In 3.1.1, we explain each step of our construction
carefully for a pseudo symplectic groupoid, and in 3.1.2, we will
deal with the Poisson case.

\subsubsection{Pseudo symplectic groupoid} Our construction of a
quantization of the noncommutative Poisson algebra is divided into
the following steps:
\begin{enumerate}
\item We use Fedosov's method to construct a resolution of a Weyl
algebra bundle and define a quantization map forgetting about the
groupoid operations, \item We prove that the quantization map is
compatible with the groupoid convolution product, \item we
construct a quantization of a groupoid algebra.
\end{enumerate}

{\bf Quantization map}

There are three steps in Fedosov's construction of a deformation
quantization of a symplectic manifold. We follow them to construct
a quantization map for a pseudo symplectic groupoid.
\begin{enumerate}
    \item A Weyl algebra bundle $\calw \index{$\calw$}$ on $\calg$;
    \item A flat connection on $\calw \index{$\calw$}$ of the form $-\delta \index{$\delta$!differential operator} +\partial +\frac {i}{h}[\ r\ ,\ \ ]$;
    \item Quantization map $\calq \index{$\calq$}: C_c ^{\infty}(\calg)\to \calw \index{$\calw$}$;
\end{enumerate}
{\bf Step 1.} We start with introducing the Weyl algebra bundle
$\calw \index{$\calw$}$ on $\calg$. We associate to every
symplectic vector space $(V, \omega)$ a Weyl algebra over
$\complex$ with a unit, which consists of formal series
\begin{equation}
\label{definition} \sum _{i\geq 0,\ |\alpha|\geq 0}\hbar ^i a_{i,
\alpha}y^{\alpha} .
\end{equation}
The multiplication $\circ \index{$\circ$}$ is defined as
\[
\begin{array}{ll}
a(\hbar, y)\circ \index{$\circ$} b(\hbar, y)&=\exp (-\frac
{i\hbar}{2}\omega ^{ij}\frac {\partial}{\partial y^i}\wedge \frac
{\partial}{\partial z^j})a(\hbar,
y)b(\hbar, z)|_{z=y}\\
&=\sum _{k=0} ^{\infty}(-\frac {i\hbar}{2})^k \frac {1}{k!}\omega
^{i_1j_1}\cdots \omega ^{i_k j_k}\frac {\partial ^k a}{\partial
y^{i_1}\cdots \partial y^{i_k}}\frac {\partial ^k b}{\partial
y^{j_1}\cdots y^{j_k}},
\end{array}
\]
where $y^i$ is the i-th coordinate on the vector space and $\omega
^{ij}$ is the inverse of the symplectic form $\omega$. It is easy
to check that the multiplication  $\circ \index{$\circ$}$ defined
above is associative and independent of choice of basis.

For a pseudo symplectic groupoid $(\calg, \calf \index{$\calf$},
\omega)$, at each point $\gamma \in \calg$, $(\calf
\index{$\calf$} |_{\gamma},\ \omega |_{\gamma})$ is a symplectic
vector space. Therefore, by the above construction, we can define
a Weyl algebra bundle by associating a Weyl algebra $\calw
\index{$\calw$}_{\gamma}$ to $(\calf \index{$\calf$} | _{\gamma},\
\omega |_{\gamma})$.

On a Weyl algebra, we prescribe a grading to the variables by
setting $deg (y^i )=1$ and $deg(\hbar) =2.$ By virtue of this
grading,  $\calw \index{$\calw$}$ becomes an $\mathbb {N}$ graded
bundle.

Here, to define a deformed groupoid algebra, we need to integrate
a Weyl algebra valued function. To make this integration
well-defined, we need to define an $\reals ^n $ translation
invariant topology on the Weyl algebra to make it a locally convex
topological space.  From (\ref{definition}), the elements $\hbar^k
y^{\alpha}$ form a basis of $\calw \index{$\calw$} _{\gamma}$.
Therefore, as a vector space, $\calw \index{$\calw$} _{\gamma}$
can be identified with $\complex ^{\mathbb {N}}$. On $\complex
^{\mathbb {N}}$, we can choose the compact open topology. The
induced compact open topology defines a topology on $\calw
\index{$\calw$} _{\gamma}$.
\begin{prop}
\label{topology}A Weyl algebra with the above topology is a
complete locally convex topological algebra.
\end{prop}
$\pf$ ``Completeness" comes from the completeness of the space
$\{f:\mathbb {N}\to \complex\}$ with the compact open topology.
The continuity of the multiplication is from the observation that
a coefficient of the multiplication is determined by a finite
number of coefficients of the two original elements. The local
convexity is a straightforward check. we know that all finite
dimensional vector spaces with the compact open topology are
locally convex spaces. The same arguments work for an infinite
dimensional vector space. $\Box$

It is easy to check that the compact open topology on a Weyl
algebra is invariant under the $GL(n, \complex)$ action on the
symplectic vector space. With the above topology of Weyl algebra,
the Weyl algebra bundle becomes a topological bundle.
\begin{rmk}
The topology we give here is not a $C^*-$norm topology. One may
define a $C^*-$norm on the Weyl algebra, but it is very hard to
make the later induction steps continuous in $C^*-$norm. The above
defined topology makes the induction continuous, but is not a
$C^*-$norm topology. There is a well-known topology---$\hbar$-adic
topology on the Weyl algebra. But it does not work in our
construction, because $\hbar$-adic topology does not make a Weyl
algebra into a locally convex topological space.
\end{rmk}
{\bf Step 2.} The next step in our construction is to find a
resolution of the Weyl algebra bundle. We show that there is a
flat connection on $\wedge ^* \calf^* \index{$\calf$} \otimes
\calw \index{$\calw$}$. For this purpose, we introduce two
operators $\delta \index{$\delta$!differential operator}$ and
$\delta \index{$\delta$!differential operator} ^*$ on $\wedge ^*
\calf^* \index{$\calf$} \otimes \calw \index{$\calw$}$, defined by
\[
\begin{array}{ll}
\delta \index{$\delta$!differential operator} a=d x^i \wedge \frac
{\partial a}{\partial y^i}, & \delta \index{$\delta$!differential
operator} ^* a=y^ki(\frac {\partial }{\partial x^k})a.
\end{array}
\]
$\delta \index{$\delta$!differential operator}\ (\delta
\index{$\delta$!differential operator} ^*)$ is a degree
decreasing\ (increasing) operator, satisfying the following
properties.
\begin{prop}
The operators $\delta \index{$\delta$!differential operator}$ and
$\delta \index{$\delta$!differential operator} ^*$ do not depend
on the choices of local coordinates and satisfy:
\begin{enumerate}
    \item $$ \delta \index{$\delta$!differential operator} ^2=(\delta \index{$\delta$!differential operator} ^*)^2=0. $$
    \item On a monomial $y^{i_1}y^{i_2}\cdots y^{i_p}dx^{j_1}\wedge dx^{j_2}\cdots
    \wedge dx^{j_q}$, we have $\delta \index{$\delta$!differential operator} \delta \index{$\delta$!differential operator} ^* +\delta \index{$\delta$!differential operator} ^* \delta \index{$\delta$!differential operator}
    =(p+q)id$. Generally, for $a\in \wedge ^* \calf \index{$\calf$}
    \otimes \calw \index{$\calw$}$,
    \[ a=\delta \index{$\delta$!differential operator} \delta \index{$\delta$!differential operator} ^{-1}a+\delta \index{$\delta$!differential operator} ^{-1}\delta \index{$\delta$!differential operator} a +a_{00}, \]
     where $\delta \index{$\delta$!differential operator} ^{-1}$ is defined by
    \[
    \begin{array}{ll}
     \delta \index{$\delta$!differential operator} ^{-1}a_{pq}
     =&\frac {1}{p+q}\delta \index{$\delta$!differential operator} ^{*}a_{pq},\
     p+q>0,\\
     \delta \index{$\delta$!differential operator} ^{-1}a_{00}=&0,
     \end{array}
    \]
    and in which $a_{pq}$ is the homogeneous part of $a$ with
    degree $p$ in $y$ and degree $q$ in $dx$.
\end{enumerate}
\end{prop}

In the following, we assume that there is an invariant pseudo
symplectic connection $\partial$ on $\calg$. For convenience, we
will omit the words ``invariant pseudo", and simply say
``symplectic connection".

A symplectic connection $\partial$ also defines a connection on
the Weyl algebra bundle and the tensor bundle $\wedge ^* \calf
\index{$\calf$} ^* \otimes \calw \index{$\calw$} |_{I}$,
\footnote{$I$ is a leave of $\calf$. For convenience, in the
following we will omit $``|_{I}"$.} which can be expressed as
\[
\partial a= dx^i \wedge \partial _i a,
\]
where $\partial _i a$ is a covariant derivative. It is not
difficult to check that $\partial$ has the following properties:
\begin{prop}
\label{prop:partial}
\begin{enumerate}
\item $\partial (a \circ \index{$\circ$} b)=\partial a\circ \index{$\circ$}  b + (-1)^{q}a \circ \index{$\circ$}
\partial b,\ \ \ \ for\ a \in \wedge ^{q}\calf \index{$\calf$} ^*\otimes \calw \index{$\calw$}$.
\item for any scalar form $\phi \in \wedge ^{q}\calf \index{$\calf$} ^{*},\
\partial(\phi \wedge a)=d\phi \wedge a + (-1)^q \phi \wedge
\partial a$.
\item $\partial \delta \index{$\delta$!differential operator} a+ \delta \index{$\delta$!differential operator} \partial a=0$. \item $\partial ^2
a=\partial (\partial a)=\frac {i}{\hbar}[R, a]$ where
\[ R=\frac {1}{4}R_{ijkl}y^iy^jdx^k\wedge dx^l,\] is the curvature
of the symplectic connection\footnote{$[\ \ ,\ \ ]$ is\ defined\
as $\  [a,b]=a\circ b-(-1)^{deg(a)deg(b)}b\circ a$. }.
\end{enumerate}
\end{prop}

In a Darboux chart, the connection can be written as
\[
\partial a=da + \frac {i}{\hbar}[\Gamma, a],
\]
where $\Gamma$ is a local 1-form with values in $\calw
\index{$\calw$}$, and $d=dx^i \wedge \frac {\partial}{\partial
x^i}$ is the exterior differential with respect to $x$. To find an
abelian connection, we will consider a connection on $\calw
\index{$\calw$}$ with a more general form,
\[
D a =\partial a +\frac {i}{\hbar}[\gamma, a]=da +\frac
{i}{\hbar}[(\Gamma +\gamma), a ],
\]
where $\gamma$ is a section of $\calw \index{$\calw$} \otimes
\wedge^1 \calf \index{$\calf$}^*$. Here, we only consider the case
where $\gamma _0 |_{y=0}=0$. For the new connection, we call the
following 2-form
\[
\Omega \index{$\Omega$!curvature} =R+ \partial \gamma +\frac
{i}{\hbar}\gamma ^2,
\]
the curvature form of $D$.
\begin{prop}
\begin{enumerate}
    \item (Bianchi identity)
    \[ D\Omega \index{$\Omega$!curvature} =\partial \Omega \index{$\Omega$!curvature} +\frac {i}{\hbar}[\gamma, \Omega \index{$\Omega$!curvature}]=0 .\]
    \item for any section $a\in \wedge ^* \calf \index{$\calf$} ^* \otimes \calw \index{$\calw$}$,
    we have
    \[ D^2 a=\frac {i}{\hbar}[\Omega \index{$\Omega$!curvature}, a] .\]
\end{enumerate}
\end{prop}
$\pf$ The proof is the same as the proof of Lemma 5.1.5 of
\cite{Fe:deformation}. $\Box$

After this long preparation, we are able to define the key notion
in Fedosov's method.
\begin{dfn}
A connection $D$ on the bundle $\calw \index{$\calw$}$ is called
abelian if for any section $a\in \wedge ^* \calf \index{$\calf$}
^* \otimes \calw \index{$\calw$}$,
\[ D^2a= \frac {i}{\hbar}[\Omega \index{$\Omega$!curvature}, a]=0 .\]
\end{dfn}
\begin{rmk}
Here, we follow Fedosov using the term ``abelian connection". This
is also known as the ``Fedosov connection" in the literature.
\end{rmk}
The main theorem in this step is
\begin{thm}
\label{conn-exist} There exists one and only one abelian
connection of the form
\[
D=-\delta \index{$\delta$!differential operator} +\partial +\frac
{i}{\hbar}[\ r,\ \ ],\] with
\[ deg (r)\geq 2,\ \ \ \ \  \delta \index{$\delta$!differential operator} ^{-1}r=0 ,\] and $r$ is
invariant under $s^*, t^*$ maps.
\end{thm}
$\pf$ The key idea of the proof is iteration. The steps are
similar to the proof of Theorem 5.2.2 in \cite{Fe:deformation}.
Here, there are two more things new from \cite{Fe:deformation}.
One is that $r$ is invariant under the $s, t$ maps. The reason for
this is that we have required all the data in the construction to
be invariant under s, t maps, i.e. the \'etalification, the 2-form
and the symplectic connection are invariant under the $s, t$ maps,
and also the iteration steps are $s, t$ invariant. Therefore, it
is straightforward to check that $r$ is also $s, t$ invariant. The
other is that all the steps of the proof in \cite{Fe:deformation}
are continuous with respect to the topology we defined in
\ref{topology}, which follows from the degree increasing iteration
procedure (see Lemma 5.2.3 in
\cite{Fe:deformation}). $\Box$\\
{\bf Step 3.} For $D$, we consider the set of flat sections $\calw
\index{$\calw$} _D\stackrel {def}{=}\{ a\in \calw \index{$\calw$}
: Da=0 \}$. It has the following important property.
\begin{prop}
\label{quantization} For any $a_0 \in C^{\infty}(\calg)[[\hbar]]$,
there exists a unique section $a \in \calw \index{$\calw$} _D$
such that $\sigma (a)=a_0$. ($\sigma (a)$ means the projection
onto the center: $\sigma(a)=a(x, 0, h)$).
\end{prop}
$\pf$ Follows from Theorem 5.2.4 of \cite{Fe:deformation}. $\Box$

By Proposition \ref{quantization}, we define a quantization map
$\calq \index{$\calq$}: C^{\infty}(\calg)[[\hbar]]\to \calw
\index{$\calw$} _D$ to be the inverse of $\sigma$ in the above
proposition.

{\bf Groupoid operation}

In the above steps, almost everything is analogous to Fedosov's
quantization of a symplectic manifold. But it only provides us
with a formal deformation quantization of the commutative algebra
$C_c^{\infty}(\calg)[[\hbar]]$ with the pointwise multiplication.
In the following, we will use this construction to quantize the
groupoid algebra $C_c^{\infty}(\calg)[[\hbar]]$ with the
convolution product. Our strategy to define the deformed groupoid
algebra is the following:
\begin{enumerate}
\item Define a new multiplication $\ast \index{$\ast$}$ on sections of Weyl
algebra bundles, which corresponds to the convolution. \item Prove
that $(C_c^{\infty}(\wedge ^* \calf \index{$\calf$} ^*\otimes
\calw \index{$\calw$}), \ast \index{$\ast$} )$ is associative and
$D$ acts as a derivation.
\end{enumerate}
\begin{rmk}
We use $\ast \index{$\ast$}$ and $\star \index{$\star$}$ for
different meanings. $\ast \index{$\ast$}$ is a new product on
smooth sections of Weyl algebra bundle, and $\star
\index{$\star$}$ is a star product on the groupoid algebra.
\end{rmk}
{\bf Step 1.} We introduce a new algebraic product $\ast
\index{$\ast$}$ on the sections of $\wedge ^*\calf \index{$\calf$}
^*\otimes \calw \index{$\calw$}$. For $f, g\in C_c ^
{\infty}(\wedge ^*\calf \index{$\calf$} ^*\otimes \calw
\index{$\calw$})$,
\[
f\ast \index{$\ast$} g(\gamma)\stackrel{def}{=}\int _{\alpha\cdot
\beta=\gamma}t^*(f(\alpha))\circ \index{$\circ$}
s^*(g(\beta))d\lambda ^{r(\gamma)}.
\]
\begin{rmk}
In the construction of $\wedge ^* \calf \index{$\calf$} ^* \otimes
\calw \index{$\calw$}$, we required everything to be invariant
under the $s,\ t$ maps. Therefore, for $\gamma \in \calg$, both
$s^* : \wedge ^* \calf \index{$\calf$} ^* \otimes\calw
\index{$\calw$} |_{s(\gamma)}\to \wedge ^* \calf \index{$\calf$}
^* \otimes \calw \index{$\calw$} | _{\gamma}$ and $t^* : \wedge ^*
\calf \index{$\calf$} ^* \otimes \calw| \index{$\calw$}
_{t(\gamma)}\to \wedge ^* \calf \index{$\calf$} ^* \otimes \calw
\index{$\calw$} |_{\gamma}$ are isomorphisms of Weyl algebras. In
this way, under $t^*$, $ f(\alpha)\in \wedge ^* \calf
\index{$\calf$} ^* \otimes \calw \index{$\calw$} |_{\alpha}$ is
mapped into $\wedge ^* \calf \index{$\calf$} ^* \otimes \calw
\index{$\calw$} |_{\gamma}$, and so is $g(\beta)$ by $s^*$. So the
$\ast \index{$\ast$}$ in the above formula is well-defined, and
when there is no confusion, we will drop $s^*, t^*$ and $\lambda
^{r(\gamma)}$.
\end{rmk}

\begin{rmk}
When integrating a Weyl algebra valued function $f$ along a
t-fiber, we have to use the compact open topology on a Weyl
algebra defined in \ref{topology}.
\end{rmk}
{\bf Step 2.} We will prove that $\ast \index{$\ast$}$ is
associative and $D$ acts as a derivation. Before doing this, we
first recall Lemma \ref{derivation}.

When restricted to the leaves of $\calf \index{$\calf$}$, $d$ has
the following formula: $\forall f,\ g\ \in C^{\infty}_c(\calg)$,
\[
d(f\diamond \index{$\diamond$} g)(\alpha)=\int _{\alpha =\beta
\gamma}(t^*)^{-1}(d f (\beta))g(\gamma)+ f(\beta)(s^*)^{-1}
d(g)(\gamma).
\]

\begin{lem}
\label{ast} $\ast \index{$\ast$}$ is an associative product on
$\Gamma(\wedge ^* \calf \index{$\calf$} ^* \otimes \calw
\index{$\calw$})$, and the connection $D$ is a derivation on
$(\Gamma(\wedge ^* \calf \index{$\calf$} ^* \otimes \calw
\index{$\calw$}), \ast \index{$\ast$})$, i.e.
\[
D(f\ast \index{$\ast$} g)=D f\ast \index{$\ast$}
g+(-1)^{deg(f)}f\ast \index{$\ast$} D g
\]
\end{lem}
$\pf$  For the associativity of $\ast \index{$\ast$}$, we compute
\[
\begin{array}{ll}
&f\ast \index{$\ast$} (g\ast \index{$\ast$} h)(\alpha)\\
=&\int _{\calg ^{s(\alpha)}}f(\alpha \beta)\circ \index{$\circ$}(g\ast \index{$\ast$} h)(\beta ^{-1})d \lambda ^{s(\alpha)}(\beta)\\
=&\int _{\calg ^{s(\alpha)}}f(\alpha \beta)d\lambda
^{s(\alpha)}(\beta) \int _{\calg ^{s(\beta ^{-1})}}(g(\beta
^{-1}\gamma)\circ \index{$\circ$} h(\gamma ^{-1}))d\lambda
^{s(\beta ^{-1})}(\gamma)\\
=&\int _{\calg ^{s(\alpha)}}f(\alpha \beta)d\lambda
^{s(\alpha)}(\beta)\int _{\calg ^{s(\alpha)}}(g(\beta
^{-1}\gamma)\circ \index{$\circ$} h(\gamma ^{-1}))d\lambda
^{s(\alpha)}(\gamma)\\
=&\int _{\calg ^{s(\alpha)}}d\lambda ^{s(\alpha)}(\beta)\int
_{\calg ^{s(\alpha)}}d\lambda ^{s(\alpha)}(\gamma)f(\alpha
\beta)\circ \index{$\circ$} (g(\beta ^{-1}\gamma)\circ \index{$\circ$} h(\gamma ^{-1}))\\
=& \int _{\calg ^{s(\alpha)}}d\lambda ^{s(\alpha)}(\gamma)\int
_{\calg ^{s(\alpha)}}d\lambda ^{s(\alpha)}(\beta)f(\alpha
\beta)\circ \index{$\circ$} (g(\beta ^{-1}\gamma)\circ \index{$\circ$} h(\gamma ^{-1}))\\
=&\int _{\calg ^{s(\alpha)}}d\lambda ^{s(\alpha)}(\gamma)\int
_{\calg ^{s(\alpha)}}d\lambda ^{s(\alpha)}(\beta)(f(\alpha
\beta)\circ \index{$\circ$} g(\beta ^{-1}\gamma))\circ \index{$\circ$} h(\gamma ^{-1})\\
=&\int _{\calg ^{s(\alpha)}}d\lambda ^{s(\alpha)}(\gamma)(f\ast
\index{$\ast$}
g)(\alpha \gamma)\circ \index{$\circ$} h(\gamma ^{-1}) \\
=&(f\ast \index{$\ast$} g)\ast \index{$\ast$} h (\alpha).
\end{array}
\]

From the formulas of $D=\partial +[r,\ ]$, to prove $D$ is a
derivation it is sufficient to show that $\partial$ and $[r,\ ]$
are derivations, respectively.
\begin{enumerate}
    \item $\partial(f\ast \index{$\ast$} g)=(\partial f)\ast \index{$\ast$} g+ (-1)^{deg(f)}f\ast \index{$\ast$} (\partial
    g)$,
    \item $[r, f\ast \index{$\ast$} g]=[r, f]\ast \index{$\ast$} g+(-1)^{deg(f)}f \ast \index{$\ast$} [r, g]$.
\end{enumerate}
By partition of unity, we write $f, g$ to be of the form
\[
f(\alpha)=\sum _{i, I, \tau}\hbar^{i}f_{i,I,\tau}(\alpha)y^{I}dx
^{\tau},\ \ \ \ \ \ g(\alpha)=\sum _{j,J, \upsilon
}\hbar^{j}g_{j,J, \upsilon}(\alpha)y^{J}dx ^{\upsilon} .
\]

For the proof of 1,
\[
\begin{array}{ll}
 &\partial(f\ast \index{$\ast$} g)(\alpha)\\
=&\partial (\int_{\beta \cdot \gamma=\alpha} f(\beta )\circ \index{$\circ$} g(\gamma)d\lambda ^{r(\alpha)} )\\
=&\partial (\int_{\beta \cdot \gamma=\alpha}\sum _{i, I,
\tau}\hbar^i f _{i, I, \tau}(\beta)y^I dx ^{\tau}\circ
\index{$\circ$} \sum _{j, J, \upsilon}\hbar^j g_{j, J,
\upsilon}(\gamma)y^J dx ^{\upsilon}d \lambda
^{r(\alpha)})\\
=&\partial(\sum _{i, I, \tau, j, J, \upsilon}\int _{\beta\cdot
\gamma=\alpha }f_{i, I, \tau}(\beta)g_{j, J,
\upsilon}(\gamma)d\lambda ^{r(\alpha)}(\hbar^iy^Idx ^{\tau})\circ
\index{$\circ$} (\hbar^j
y^J dx ^{\upsilon}))\\
=&\partial(\sum _{i, I, \tau, j, J, \upsilon}f_{i, I,
\tau}\diamond \index{$\diamond$} g_{j, J, \upsilon}
(\alpha)(\hbar^iy^Idx ^{\tau})\circ \index{$\circ$}
(\hbar^j y^J dx ^{\upsilon}))\\
=&\sum _{i, I, \tau, j, J, \upsilon}\partial (f_{i, I,
\tau}\diamond \index{$\diamond$} g_{j, J, \upsilon}
(\alpha)(\hbar^iy^Idx ^{\tau})\circ \index{$\circ$}
(\hbar^j y^J dx ^{\upsilon}))\\
=&\sum _{i, I, \tau, j, J, \upsilon}(\partial (f_{i, I,
\tau}\diamond \index{$\diamond$} g_{j, J,
\upsilon})(\alpha)(\hbar^iy^Idx ^{\tau})\circ \index{$\circ$}
(\hbar^j y^J dx ^{\upsilon})\\
+&f_{i, I, \tau}\diamond \index{$\diamond$} g_{j, J,
\upsilon}(\alpha)\partial ((\hbar^iy^Idx ^{\tau})\circ
\index{$\circ$} (\hbar^j y^J dx
^{\upsilon}))),\\
\end{array}
\]
By Lemma \ref{derivation}, and the derivation property of
$\partial$ in Proposition \ref{prop:partial},
\[
\begin{array}{ll}
=&\sum _{i, I, \tau, j, J, \upsilon} (\partial(f_{i, I,
\tau})\diamond \index{$\diamond$} g_{j, J,
\upsilon})(\alpha)+f_{i, I, \tau}\diamond \index{$\diamond$}
(\partial {g_{j, J, \upsilon}})(\alpha))(\hbar^iy^Idx
^{\tau})\circ \index{$\circ$}
(\hbar^j y^J dx ^{\upsilon})\\
+&f_{i, I, \tau}\diamond \index{$\diamond$} g_{j, J,
\upsilon}(\alpha)((\partial(\hbar^iy^Idx ^{\tau}))\circ
\index{$\circ$} (\hbar^j y^J dx
^{\upsilon})+(-1)^{\tau}(\hbar^iy^Idx ^{\tau})\circ
\index{$\circ$} (\partial(\hbar^j y^J
dx^{\upsilon})))\\
=&\sum _{i, I, \tau, j, J, \upsilon} (\partial(f_{i, I,
\tau})\diamond \index{$\diamond$} g_{j, J,
\upsilon})(\alpha)(\hbar^iy^Idx ^{\tau})\circ \index{$\circ$}
(\hbar^j y^J dx ^{\upsilon})\\
+&f_{i, I, \tau}\diamond \index{$\diamond$} g_{j, J,
\upsilon}(\alpha)(\partial(\hbar^iy^Idx ^{\tau}))\circ
\index{$\circ$} (\hbar^j y^J dx ^{\upsilon})\ +f_{i, I,
\tau}\diamond \index{$\diamond$} (\partial {g_{j, J,
\upsilon}})(\alpha)(\hbar^iy^Idx ^{\tau})\circ \index{$\circ$}
(\hbar^j y^J dx
^{\upsilon})\\
+&(-1)^{\tau}f_{i, I, \tau}\diamond \index{$\diamond$} g_{j, J,
\upsilon}(\alpha) (\hbar^iy^Idx ^{\tau})\circ \index{$\circ$}
(\partial(\hbar^j y^J dx
^{\upsilon}))\\
=&(\partial f)\ast \index{$\ast$} g+ (-1)^{deg(f)}f\ast
\index{$\ast$} (\partial
    g).
\end{array}
\]
The proof for $[r, \ ]$ is much easier than the proof for $D$,
which can be derived from the invariance of $r$ under $s^*,
t^*$(see Theorem \ref{conn-exist}) and associativity of $\circ
\index{$\circ$}$.\ \ \ $\Box$

{\bf Quantization of a groupoid algebra}

From Lemma \ref{ast}, we know that sections of $\calw
\index{$\calw$} _{D}$ are closed under convolution $\ast
\index{$\ast$}$, and form a new associative algebra. Therefore, we
can define a formal deformation quantization of the groupoid
algebra $C_{c}^{\infty}(\calg)[[\hbar]]$ to be
\begin{equation}
\label{mult}
 f \star \index{$\star$} g(\gamma)=\sigma (\int _{\alpha
\cdot \beta =\gamma} t^*(\calq \index{$\calq$} (f)(\alpha))\circ
\index{$\circ$} s^* (\calq \index{$\calq$} (g)(\beta))d\lambda
^{t(\gamma)} ).
\end{equation}
which can also be written as $\sigma(\calq \index{$\calq$} (f)\ast
\index{$\ast$} \calq \index{$\calq$} (g))$.
\begin{thm}
\label{thm:star-sympl} $``\star \index{$\star$} "$ is associative.
\end{thm}
$\pf$ The associativity is an easy corollary from the
associativity of $\ast \index{$\ast$}$, and the fact that
\[
D(\calq \index{$\calq$} (f)\ast \index{$\ast$} \calq
\index{$\calq$} (g))=D(\calq \index{$\calq$} (f))\ast
\index{$\ast$} \calq \index{$\calq$} (g)+ \calq \index{$\calq$}
(f)\ast \index{$\ast$} D(\calq \index{$\calq$} (g))=0.\ \ \ \ \
\Box .
\]
Now we have constructed an associative star product on the
groupoid algebra. To show that it is a deformation quantization,
we have to prove that the first term of its semiclassical
expansion is the noncommutative Poisson structure defined in
(\ref{Poisson}). The steps to prove this are the same as in the
Poisson case, so we leave it to the next subsection.

\begin{rmk}
One can also consider the classification of the constructed
algebra, either up to isomorphism as in \cite{Fe:deformation} or
Morita equivalence \index{Morita equivalence} as in
\cite{b1:morita}.

Here, we follow the definition of a formal deformation
quantization in \cite{Fe:deformation}, in which there is no
concern about the involution on $C_{c}^{\infty}(\calg)[[\hbar]]$.
One can follow Neumaier's work \cite{n:involution} to define an
involution on the deformed algebra. Furthermore, one can also look
for a  positive (Hermitian) deformation quantization of a groupoid
algebra. We leave all this for future study.
\end{rmk}

\subsubsection{Pseudo Poisson groupoid} We have already seen a
deformation quantization of a pseudo symplectic groupoid. In this
subsection, we will work on a pseudo Poisson groupoid. The spirit
of the construction is similar to the symplectic case. However,
there are two major differences. One is that we will replace the
Weyl algebra by the Kontsevich algebra constructed by Kontsevich
in \cite{Ko:deformation} for a Poisson structure on $\reals ^n$,
the other is that on a general pseudo Poisson groupoid there is no
Poisson connection. To solve this difficulty, we will use a
structure other than a connection
---a ``quasi-connection". To lift this to the Kontsevich
algebra bundle, we will apply Kontsevich's local formality
\index{formality} theorem. Since Theorem \ref{thm:coboundary} has
already assumed the condition that $\calg$ is proper, we will
assume in this subsection that $\calg$ is proper, which in the
case of pseudo symplectic (or regular Poisson) groupoids is
stronger than the existence of an invariant connection. The idea
of using Fedosov's method to construct a deformation quantization
for a general Poisson manifold comes from Cattaneo, Felder and
Tomassini in \cite{Cf1:local-global}. As most of the steps are
similar to those in the symplectic case, we will be a little bit
sketchy. Readers are referred to \cite{Cf1:local-global} and
\cite{Ko:deformation} for details.

Essentially, there are four steps in our construction of a formal
deformation quantization of a pseudo Poisson groupoid.
\begin{enumerate}
    \item the Kontsevich algebra $\frakk \index{$\frakk$}$;
    \item the quantization bundle;
    \item the abelian connection;
    \item the quantization map $\calq \index{$\calq$}$.
\end{enumerate}
{\bf Step 1.} Let us start by reviewing Kontsevich's star product
and formality \index{formality} theorem on $\reals ^n$.
Kontsevich's formality \index{formality} theorem constructs a
quasi-isomorphism between the DGLA of the Hochschild complex of
the algebra of polynomials and the DGLA of multi-vector fields.
Kontsevich defines  a sequence of operators $U_j$ to be a
multi-linear symmetric function of $j$ arguments of multi-vector
fields, with values in the multi-differential operators
$C^{\infty}(\reals ^n)^{\otimes r}\to C^{\infty}(\reals ^n)$,
where $r=\sum _k m_k-2j+2$. The maps are $GL(d, \reals)$-invariant
and satisfy the following famous equalities:
\begin{thm}
\label{formality} (\cite{Ko:deformation}, Theorem 3.1 of
\cite{Cf1:local-global}) Let $\alpha _j \in \Gamma (\reals ^n,
\wedge ^{m_j}T \reals ^n),\ j=1, \cdots , s$  be multi-vector
fields.

Let $\epsilon _{ij}=(-1)^{(m_1+\cdots +m_{i-1})m_i+(m_1+\cdots +
m_{i-1} +m_{i+1}+\cdots +m_{j-1})m_j}$. Then, for any functions
$f_0, \cdots, f_m ,$
\[
\begin{array}{ll}
\sum _{l=0} ^{s} &\sum _{k=-1} ^{m} \sum
_{i=0}^{m-k}(-1)^{k(i+1)+m}\sum _{\sigma \in S_{i, s-l}}\epsilon
(\sigma )U_{l}(\alpha _{\sigma (1)}, \cdots, \alpha _{\sigma
(l)})(f_0 \otimes \cdots \\ &\otimes f_{i-1}\otimes U_{s-l}(\alpha
_{\sigma (l+1)}, \cdots , \alpha _{\sigma (s)})(f_i \otimes \cdots
\otimes f_{i+k})\otimes f_{i+k+1}\otimes
\cdots \otimes f_{m})\\
&=\sum _{i<j}\epsilon _{ij}U_{s-1}([ \alpha _i, \alpha _j], \alpha
_1, \cdots, \hat {\alpha _i}, \cdots , \hat {\alpha _j},\cdots,
\alpha _s)(f_0\otimes \cdots \otimes f_m).
\end{array}
\]
\end{thm}
To formally deformation quantize a Poisson manifold and a pseudo
Poisson groupoid, we will only use some special cases of the above
theorem. We follow the notation in \cite{Cf1:local-global}.
\[
\begin{array}{ll}
P(\pi) \index{$P(\pi)$}&= \sum _{j=0}^{\infty}\frac
{\hbar^j}{j!}U_j (\pi, \cdots ,
\pi)\\
A(\xi , \pi) \index{$A(\xi , \pi)$}&= \sum _{j=0} ^{\infty}\frac
{\hbar^j}{j!}U_{j+1}(\xi,
\pi, \cdots, \pi)\\
F(\xi, \eta, \pi) \index{$F(\xi, \eta, \pi)$}&= \sum
_{j=0}^{\infty}\frac {\hbar^j}{j!}U_{j+2}(\xi, \eta, \pi, \cdots,
\pi).
\end{array}
\]
$\xi , \eta$ are vector fields on $\reals ^n$, and $\pi$ is a
Poisson bivector.

Let $\mathfrak {U} \index{$\mathfrak {U}$}$ be the set of
polynomial maps from $\Gamma (\reals ^n, \wedge ^2 T\reals ^n)$ to
the space of multi-differential operators on $\reals ^n$,  with
differential $\delta \index{$\delta$!differential operator}$
defined by:
\[
\begin{array}{ll}
\delta \index{$\delta$!differential operator} S(\xi _1, \cdots,
\xi _{p+1}, \pi)=&-\sum _{i=1}^{p+1}\frac {d}{dt}|_{t=0}S(\xi _1,
\cdots, \hat {\xi _i}, \cdots, \xi _{p+1},
(\Phi _{\xi _i}^t)_* \pi)\\
&+\sum _{i<j}(-1)^{i+j}S([\xi _i, \xi _j], \xi _1, \cdots , \hat
{\xi _{i}}, \cdots, \hat {\xi _{j}}, \cdots, \xi _{p+1}, \pi),
\end{array}
\]
where $(\Phi ^{t}_{\xi _i})$ stands for the flow generated by the
vector field $\xi _i$.

The formality theorem (Theorem \ref{formality}) implies that
$P(\pi)$, $A(\xi, \pi)$, and $F(\xi, \eta, \pi)$ satisfy the
following equations.
\begin{lem}
\label{lem:curvature} (Corollary 3.2, \cite{Cf1:local-global} )
Theorem \ref{formality} implies the following:
\[
\begin{array}{l}
(i)P(\pi) \index{$P(\pi)$}\cdot (A(\xi, \pi)\otimes Id +Id \otimes
A(\xi,
\pi))-A(\xi, \pi)\cdot P(\pi) \index{$P(\pi)$}=\delta \index{$\delta$!differential operator} P(\xi, \pi).\\
(ii)P(\pi) \index{$P(\pi)$}\cdot(F(\xi, \eta , \pi)\otimes Id - Id
\otimes F(\xi, \eta, \pi) \index{$F(\xi, \eta, \pi)$})-A(\xi ,
\pi)\cdot A(\eta, \pi)
+A(\eta, \pi)\cdot A(\xi, \pi)=\delta \index{$\delta$!differential operator} A (\xi, \eta , \pi).\\
(iii)-A(\xi, \pi)\cdot F(\eta, \zeta, \pi)-A(\eta, \pi)\cdot
F(\zeta, \xi, \pi)-A(\zeta, \pi)\cdot F(\xi, \eta, \pi  ) =\delta
\index{$\delta$!differential operator} F(\xi, \eta, \zeta, \pi).
\end{array}
\]

With Lemma \ref{lem:curvature}, we define a $\circ \index{$\circ$}
$ product on $(\reals ^n, \pi)$,
\[
f\circ \index{$\circ$} g\stackrel {def}{=}P(\pi)
\index{$P(\pi)$}(f\otimes g),\ \ \ \ f, g \in \reals [[y^1,
\cdots, y^n]][[\hbar]],
\]
which is associative. We call $\reals [[y^1, \cdots,
y^n]][[\hbar]]$, with the $\circ \index{$\circ$}$ product the
Kontsevich algebra $\frakk \index{$\frakk$}$.
\end{lem}
\begin{rmk}
Here, we work with $\reals-$algebras. One can generalize this to
$\complex-$algebras without any extra effort. Also, one can define
a topology as in the last section.
\end{rmk}
{\bf Step 2.} Parallel to the construction of the Weyl algebra
bundle, we define a Kontsevich algebra bundle. The construction is
a little more involved than the Weyl algebra bundle, because the
Kontsevich algebras at different points are generally not
isomorphic. Therefore, we must use some formal geometry.

Let $(\calg, \calf \index{$\calf$}, \pi)$ be a proper pseudo
Poisson groupoid, and at each point $\gamma$ of $\calg$, let $L
_{\gamma}$ be the leaf of $\calf \index{$\calf$}$ through $\gamma$
with rank $n$. With the help of formal geometry, we can define a
vector bundle $E$ on $\calg$, with fiber $\reals [[ y^1, \cdots,
y^n]][[\hbar]]$(\cite{Cf1:local-global}). By the same arguments as
in \cite{Cf1:local-global}, namely the triviality of an infinite
jet bundle, one can find a family $(\phi _{\gamma})_{\gamma \in
\calg}$ of infinite jets at the zero of local diffeomorphisms
$\phi _{\gamma}: (\reals ^n, 0) \to (L_{\gamma}, \gamma)$ such
that $\phi _{\gamma}(0)=\gamma$. This family of infinite jets at
zero is called a {\bf quasi-connection}. Because $\phi _{\gamma}$
is a local diffeomorphism, it pulls back the Poisson structure
$\pi$ to $\reals ^n$. Therefore, at fiber $E_{\gamma}$, using
Kontsevich's product, one can define a Kontsevich algebra $\frakk
\index{$\frakk$} _{\gamma}$. In this way, we obtain a Kontsevich
algebra bundle $\frakk \index{$\frakk$}$. Here, because $\calg$ is
proper, we can integrate $\phi _{\gamma}$ along
$\calg$ to get a $s, t$ invariant quasi-connection.\\
{\bf Step 3. } As in the pseudo symplectic groupoid case, we
consider $ \wedge ^* \calf \index{$\calf$} ^* \otimes \frakk
\index{$\frakk$}$ in which $\wedge ^* \calf \index{$\calf$} ^*$
denotes differential forms along $\calf \index{$\calf$}$. We
introduce a connection $D: \Gamma (\frakk \index{$\frakk$})\to
\Gamma (\frakk \index{$\frakk$})$ by
\[
(Df)_{\gamma}=d_x f + A_x ^{M}f,
\]
where the $d_x f$ is the de Rham differential of $f$, viewed as a
function on $\calg$, and $A_x ^{M}$ is defined in
\cite{Cf1:local-global}.  From the invariance of $\phi$ under the
$s, t$ maps, $D$ is also invariant under $s, t$. Furthermore, we
have the following proposition for $D$,
\begin{prop}
Let $F^{\calg}$ be a $\frakk \index{$\frakk$}$ valued two form on
$E$, defined by
\[
F^{\calg}=F((\hat{\xi}_{\gamma}),(\hat{\eta} _{\gamma}), \pi
_{\gamma}).
\]
where $\hat{\xi}_{\gamma}$ and $\hat{\eta}_{\gamma}$ are lifted
vectors on $E$ as in \cite{Cf1:local-global}. Then for any $f,
g\in \Gamma (\frakk \index{$\frakk$})$,
\begin{enumerate}
\item $ D(f\circ \index{$\circ$} g)=Df\circ \index{$\circ$} g + f\circ \index{$\circ$} D g, $ \item
$D^2f=F^{\calg}\circ \index{$\circ$} f-f \circ \index{$\circ$}
F^{\calg}, $ \item $ D F^{\calg}=0 .$
\end{enumerate}
\end{prop}
\begin{rmk}Here, the $\circ \index{$\circ$}$ between $F^{\calg}$ and $f$ is the
fiberwise multiplication on $\wedge ^* \calf \index{$\calf$} ^*
\otimes E$. It is easy to check that the multiplication is well
defined.
\end{rmk}
$\pf$ Same as the proof of Proposition 4.2 in
\cite{Cf1:local-global}.

Using $D$, with the same arguments of \cite{Cf1:local-global}, we
have the following theorem:
\begin{thm}
There exists an abelian connection $\bar {D}$ of the form $D+[\
r,\ \ ]=D_0 + \hbar D_1 +\hbar^2 D_2 +\cdots $ on $\frakk
\index{$\frakk$}$, so that there is an isomorphism $\calq
\index{$\calq$} : C_{c}^{\infty}(\calg)[[\hbar]] \to H^0(\frakk
\index{$\frakk$}, \bar {D})$.
\end{thm}
$\pf$ The proof follows \cite{Cf1:local-global}.
\begin{rmk}
Here,  by integration along $\calg$, the construction of $r$ and
$\calq \index{$\calq$}$ can be made invariant under the maps $s,
t$.
\end{rmk}
{\bf Step 4.} On $\frakk \index{$\frakk$}$, we introduce a new
algebraic structure, for $f, g\in \frakk \index{$\frakk$}$,
\[
f\ast \index{$\star$} g(\alpha)\stackrel {def}{=}\int_{\beta \cdot
\gamma =\alpha} t^*(f(\beta)) \circ \index{$\circ$} s^*
(g(\gamma)) .
\]
Similar to \ref{ast}, we have the following properties of $\ast
\index{$\ast$}$.
\begin{prop}
\label{prop:star-poi-conn}
\begin{enumerate}
\item $(f\ast \index{$\ast$} g)\ast \index{$\ast$} h=f\ast \index{$\ast$} (g\ast \index{$\ast$} h),$ \item $D(f\ast \index{$\ast$} g)=Df
\ast \index{$\ast$} g+ f\ast \index{$\ast$} Dg,$ \item $[r, f\ast
\index{$\ast$} g]=[r, f]\ast \index{$\ast$} g+ (-1)^{deg(f)}f\ast
\index{$\ast$} [r, g]$.
\end{enumerate}
\end{prop}

From Proposition \ref{prop:star-poi-conn}, we know that the kernel
of $\bar {D}$ is closed under $\ast \index{$\ast$} $. Therefore,
the deformed groupoid algebra $\calg$ is defined as
$(C_{c}^{\infty}(\calg)[[\hbar]], \star \index{$\star$})$, for $f,
g\in C_{c}^{\infty}(\calg)[[\hbar]]$,
\[
f\star \index{$\star$} g(\alpha)\stackrel {def}{=}\calq
\index{$\calq$}^{-1}(\int_{\beta \cdot \gamma =\alpha} t^*(\calq
\index{$\calq$} (f)(\beta)) \circ \index{$\circ$} s^* (\calq
\index{$\calq$} (g)(\gamma))) .
\]

It is easy to check that $\star \index{$\star$} $ is associative
from the first equality of Proposition \ref{prop:star-poi-conn}.
To finish the proof that $\star \index{$\star$}$ defines a
deformation quantization of a groupoid algebra, we still need to
show that the linearization of the $\star \index{$\star$}$ is the
noncommutative Poisson structure $\Pi$ defined in Section 2.

As $\calq \index{$\calq$} ^{-1}$ is evaluation at $y=0$, which is
independent of integration along $s$ and $t$ fibers,
\[
\begin{array}{ll}
&f\star \index{$\star$} g(\alpha) \\
=& \int _{\alpha =\beta \gamma}\calq \index{$\calq$}
^{-1}(t^*(\calq \index{$\calq$} (f)(\beta))\circ \index{$\circ$}
s^*(\calq \index{$\calq$} (g)(\gamma))).
\end{array}
\]
It is easy to check that $\calq \index{$\calq$}$ is invariant
under the $s,\ t$ maps, so
\[
f\star \index{$\star$} g(\alpha)=\int _{\alpha =\beta \gamma}\calq
\index{$\calq$} ^{-1}(\calq \index{$\calq$} (t^* (f(\beta)))\circ
\index{$\circ$} \calq \index{$\calq$} (s^* (g(\gamma)))).
\]
By this expression, and the fact that $\circ \index{$\circ$}$ is a
formal deformation quantization of the pointwise  multiplication,
we have
\[
\begin{array}{ll}
&f\star \index{$\star$} g(\alpha) \\
=&\int _{\alpha =\beta \gamma} i\hbar\pi (t^*(f(\gamma)), s^*
(g(\beta)))+o(\hbar)\\
=&i\hbar \Pi(f, g)(\alpha)+o(\hbar)\ \ \ \ \ \ \Box
\end{array}
\]

In summary, we have proved the following theorem.
\begin{thm}
For a proper pseudo Poisson groupoid, there is always a formal
deformation quantization of the noncommutative Poisson algebra
defined in Section 2.
\end{thm}

\begin{cor}
A Poisson orbifold groupoid can be formally deformation quantized
by the methods used in this section.
\end{cor}
\subsection{ A trace formula} Traces on an algebra play an
important role in understanding the algebra. The number of traces
and the value of a trace provide a large amount of information
about an algebra. In the study of the quantum index theorem, we
will consider the value of a trace on a quantized algebra valued
projection matrix, which requires the study of properties of
traces on a quantized groupoid algebra. In this section, we will
show a construction which defines a natural trace on a deformation
quantized groupoid algebra. In \cite{ta:trace}, we will study the
Hochschild and cyclic homology of a quantized groupoid algebra
which consists of generalized traces.

A trace of a formal algebra $A[[\hbar]]$ is defined as a $\complex
[[\hbar ]]\ (\reals [[\hbar ]])$ valued linear functional on
$A[[\hbar ]]$, vanishing on commutators. Traces on a deformation
quantization of a symplectic (Poisson) manifold have been well
studied. It is shown that on a quantized symplectic manifold there
is a unique trace up to normalization. Furthermore, following an
idea from noncommutative geometry, Connes, Flato and Sternheimer
introduced the notion of a closed $\star \index{$\star$}-$product.
The existence of a closed $\star \index{$\star$}-$product on a
symplectic manifold was first proved by Omori, Maeda and Yoshioka,
and on a Poisson manifold the corresponding result is due to the
work of Felder and Shoikhet in \cite{FS:trace-Poi}.

In the following, we define and prove the existence of a closed
deformation quantization of a pseudo symplectic (Poisson) groupoid
\index{pseudo symplectic (Poisson) groupoid} with an invariant
measure.
\begin{dfn}
A formal deformation quantization $(C_c ^{\infty}(\calg)[[\hbar]],
\star \index{$\star$})$ of a pseudo symplectic (Poisson) groupoid
\index{pseudo symplectic (Poisson) groupoid} $\calg$ is closed if
there is a groupoid invariant volume form\footnote{Usually, if a
groupoid has a symmetric measure, then the same formula defines a
trace on $C_c^{\infty}(\calg)$. Here, we require a stronger
condition --``invariant measure" for later use.} $\Omega
\index{$\Omega$!volume}$ on $\calg ^{(0)}$, such that
\[
Tr_{\Omega \index{$\Omega$!volume}}\stackrel{def}{=}\int _{\calg
^{(0)}} f |_{\calg ^{(0)}}\Omega \index{$\Omega$!volume}
\]
is a trace on the deformed algebra\footnote{It is easy to check
that the same formula also defines a trace on the groupoid
algebra. }.
\end{dfn}

The main result of this section is the following theorem.
\begin{thm}
\label{thm:trace}Let $\calg$ be a pseudo symplectic (Poisson)
groupoid \index{pseudo symplectic (Poisson) groupoid} with an
invariant connection\footnote{If $\calg$ is a general pseudo
Poisson groupoid, then instead of the existence of an invariant
connection, we require $\calg$ to be proper.}, and $\Omega
\index{$\Omega$!volume}$ an invariant volume form on the unit
space $\calg^{(0)}$, satisfying $div_{\Omega
\index{$\Omega$!volume}}\pi =0$. Then
\begin{equation}
\label{trace}  \int _{\calg ^{(0)}}f \Omega
\index{$\Omega$!volume}
\end{equation}
defines a trace on the deformation quantized groupoid algebra
$(C_c^{\infty}(\calg)[[\hbar]], \star \index{$\star$})$.
\end{thm}
\begin{rmk}
The $\star \index{$\star$}$ product might be different from the
original construction in Section 3.1.
\end{rmk}
\begin{rmk}
For the symplectic case, one has a natural choice of $\Omega
\index{$\Omega$!volume}$ which is $\omega ^{\frac {1}{2}dim(\calg
^{(0)})}$. To prove that it defines a trace, one can first follow
the method in \cite{Fe:deformation} to find a closed star product.
Then the following arguments show that $\Omega
\index{$\Omega$!volume} ^{\frac {1}{2}dim(\calg _0)}$ defines a
trace.
\end{rmk}
\begin{rmk}
For a pseudo Poisson groupoid, we can define three ``modular
classes". One is the modular class of the foliation generated by
$\calf \index{$\calf$}$, which is 0 if and only if there exists a
transversal Haar system;\index{transversal Haar system}
\index{Haar system} one is the modular class of the groupoid,
which vanishes if and only if there is an invariant Haar system;
\index{Haar system} the other is the modular class of the Poisson
structure, which is zero if and only if there is a volume $\Omega
\index{$\Omega$!volume}$ on $\calg^{(0)}$ with $div_{\Omega
\index{$\Omega$!volume}}\pi =0$. Therefore, in Theorem
\ref{thm:trace}, we are working with a groupoid whose all three
modular classes are trivial.
\end{rmk}
\begin{rmk}
We will show in \cite{ta:trace} that a quantized pseudo symplectic
(Poisson) groupoid \index{pseudo symplectic (Poisson) groupoid}
usually has more than one trace, the number of which is determined
by the 0-th Hochschild homology of the quantized groupoid. Formula
(\ref{trace}) provides one example of traces. It would be
interesting to find formulas of other traces on a deformation
quantized pseudo symplectic (Poisson) groupoid \index{pseudo
symplectic (Poisson) groupoid}.
\end{rmk}
$\pf$ We will work on the Poisson case.

As $\Omega \index{$\Omega$!volume}$ is invariant under groupoid
action, we may define an $n-$form $s^* \Omega
\index{$\Omega$!volume}(=t^*\Omega \index{$\Omega$!volume})$ on
$\calg$, $n=dim(\calg^{(0)})$. Because of the invariance
assumption, it is easy to check that on the integrated submanifold
$L_{\gamma}$, the divergence of $\pi$ with respect to $s^*\Omega
\index{$\Omega$!volume}$ is still $0$. For $\Omega
\index{$\Omega$!volume}$, following the same method of
\cite{FS:trace-Poi}, we can construct a family of
$L_{\infty}-$morphisms
\[
U_{L_{\Gamma}}: [T^* _{poly}(L_{\Gamma})]_{div}\to [D^*
_{poly}(L_{\gamma})]_{cycl},
\]
where $[T^* _{poly}(L_{\Gamma})]_{div}$ is the subalgebra of the
DGLA of multi-vector fields whose elements have 0 divergence and
$[D^* _{poly}(L_{\gamma})]_{cycl}$ is the DGLA  of cyclic cochain
complex (see \cite{FS:trace-Poi}). From this map, a divergence 0
Poisson structure defines a product $\circ \index{$\circ$}$ on
$C_c^{\infty}(L_{\gamma})[[\hbar]]$ such that for any three smooth
functions $f, g, h$ with compact supports, we have:
\begin{equation}
\label{eq:cyclic} \int _{L_{\gamma}}(f\circ \index{$\circ$}
g)h\Omega \index{$\Omega$!volume}=\int _{L_{\gamma}}(g\circ
\index{$\circ$} h)f\Omega \index{$\Omega$!volume}.
\end{equation}

Using $\circ \index{$\circ$}$ and following the same arguments as
in the last section, we can easily prove that
\[
f\star \index{$\star$} g(\alpha)=\calq \index{$\calq$} ^{-1}(\int
_{\beta \cdot \gamma =\alpha}t^*(\calq \index{$\calq$}
(f)(\beta))\circ \index{$\circ$} s^*(\calq
\index{$\calq$}(g)(\gamma))d\lambda ^{r(\alpha)})
\]
defines a star product on the groupoid algebra
$C_c^{\infty}(\calg)[[\hbar]]$. In the following, we prove that
(\ref{trace}) defines a trace on this star product.

First, we know that $\calq \index{$\calq$} ^{-1}$ is evaluation at
$y=0$, so there is no difference if we put $\calq \index{$\calq$}
^{-1}$ inside the integral. Then $Tr_{\Omega}(f\star g)=\int
_{\calg _0}f\star \index{$\star$} g \Omega
\index{$\Omega$!volume}$ can be written as
\[
\int _{\calg _0}\Omega \index{$\Omega$!volume} \int _{\alpha \cdot
\alpha ^{-1}=x}\calq \index{$\calq$} ^{-1}(t^*(\calq
\index{$\calq$} (f)(\alpha))\circ \index{$\circ$} s^*(\calq
\index{$\calq$} (g)(\alpha ^{-1})))d\lambda ^{x}.
\]
If we set $\tilde {g}(\alpha)\stackrel {def}{=}g(\alpha ^{-1})$,
by the invariance of $\calq \index{$\calq$}$ under groupoid
operations, e.g. $s^*,\ t^*$, and inverse maps, the trace can
again be written using $\circ \index{$\circ$}$,
\[
\int _{\calg _0}\Omega \index{$\Omega$!volume} \int _{\alpha \cdot
\alpha ^{-1}=x}t^*(f\circ \index{$\circ$} \tilde
{g}(\alpha))d\lambda ^{x}. \footnote{Here, it is safe to discard
$t^*$, because $f\ast g$ is a formal function and not a section.}
\]

Since $\calf \index{$\calf$}$ is transverse to the t-fiber,
$t^*(\Omega \index{$\Omega$!volume})$ and $\lambda ^{x}$ together
forms a Borel measure $\tilde {\Omega \index{$\Omega$!volume}}$ on
$\calg$. The above can be summarized as
\[
\int _{\calg _0}f\star \index{$\star$} g \Omega
\index{$\Omega$!volume} =\int f\circ \index{$\circ$} \tilde {g}
\tilde {\Omega \index{$\Omega$!volume}} .
\]

On the other hand, according to the invariance of the Haar system
\index{Haar system}, one can choose a transversal
sub-manifold\footnote{Here, one might need to work in local charts
and use partition of unity.} $V$ to the foliation of $\calf
\index{$\calf$}$, and by Fubini's theorem, the integral $\int
f\circ \index{$\circ$} g \tilde {\Omega \index{$\Omega$!volume}}$
can be written as
\[
\int _{V} \lambda _{V} \int _{L_{\gamma}}f\circ \index{$\circ$}
\tilde {g} \Omega \index{$\Omega$!volume} .
\]
By setting $h=1$ in equation (\ref{eq:cyclic}), we have
\[
\int _{L_{\gamma}}f\circ \index{$\circ$} \tilde {g} \Omega
\index{$\Omega$!volume} =\int _{L_{\gamma}}f \tilde {g} \Omega
\index{$\Omega$!volume} .
\]
Therefore,
\[
\int f\circ \index{$\circ$} \tilde {g} \tilde {\Omega
\index{$\Omega$!volume}}=\int _{V} \lambda _{V} \int
_{L_{\gamma}}f \tilde {g} \Omega \index{$\Omega$!volume},
\]
which is equal to
\[
\int _{\calg _0}\Omega \index{$\Omega$!volume} \int _{\alpha \cdot
\alpha ^{-1}=x}f(\alpha)g(\alpha ^{-1})d\lambda ^{x}=Tr_{\Omega
\index{$\Omega$!volume}}(f\diamond \index{$\diamond$} g).
\]

The last line is our definition of $Tr_{\Omega
\index{$\Omega$!volume}}(f\diamond \index{$\diamond$} g)$, which
is the same as $Tr_{\Omega \index{$\Omega$!volume}}(g\diamond
\index{$\diamond$} f)$. Therefore, what we have shown is that
$Tr_{\Omega \index{$\Omega$!volume}}$ is also a trace on the
deformed groupoid algebra. $\Box$

\subsection{Example: transformation groupoids} Let $\calf
\index{$\calf$}$ be the $ M$ component of the tangent bundle of $T
(M \times G)$. It is not difficult to check that $\calf
\index{$\calf$}$ forms an \'etalification \index{\'etalification}
on $M\rtimes G$. At each $t-$fiber, we fix the Haar measure
$\lambda$ of $G$, which is constant along $\calf \index{$\calf$}$.
It is easy to check that $\pi$ on $\calf \index{$\calf$}$ makes
the transformation groupoid $M\rtimes G\rightrightarrows M$ into a
pseudo Poisson groupoid, which induces a noncommutative Poisson
structure on its groupoid algebra $C_c ^{\infty} (G, C_c
^{\infty}(M))$, i.e. $\forall f, g \in C_c (G, C_c ^{\infty}(M))$
\[
\Pi(f, g)(x, \alpha)=\int_G  \{ f(\beta \cdot x, \alpha \beta
^{-1}), g(x, \beta)\}d\lambda(\beta).
\]

If the $G$ action is proper, $M\rtimes G\rightrightarrows M$ is
also proper, and has a quasi-connection. Therefore, $M\rtimes
G\rightrightarrows M$ can be formally deformation quantized.
\[
f\star \index{$\star$} g (m, \gamma)=\calq \index{$\calq$}
^{-1}\{\int _{\alpha \cdot \beta=\gamma} t^*(\calq \index{$\calq$}
(f)(\beta \cdot m, \alpha))\circ \index{$\circ$} s^*(\calq
\index{$\calq$} (g)(m, \beta))\lambda \}.
\]

One can directly check the following facts about the above
quantization, which we will not prove.
\begin{enumerate}
\item The nontrivial part of the above star product is on the $M$
components. \item The construction of the star product is
invariant under the $G$ action. \item $\calq \index{$\calq$}
^{-1}$ commutes with integration along $G$.
\end{enumerate}
Using the above facts, we can rewrite the star product as
\[
\begin{array}{ll}
&f\star \index{$\star$} g(m, \gamma)\\
=&\int _{\alpha\cdot \beta =\gamma} \calq \index{$\calq$}
^{-1}(t^*(\calq \index{$\calq$}
(f)(\beta \cdot m, \alpha))\circ \index{$\circ$} s^*(\calq \index{$\calq$} (g)(m, \beta)))\lambda \\
=&\int _{\alpha\cdot \beta =\gamma} \calq \index{$\calq$}
^{-1}(\calq \index{$\calq$} (\beta
^*(f(\alpha)))(m)\circ \index{$\circ$} \calq \index{$\calq$} (g(\beta))(m))\lambda\\
=&\int _{\alpha\cdot \beta =\gamma} \beta ^*(f(\alpha))\star
\index{$\star$} g(\beta)(m)\lambda .
\end{array}
\]
In the second step, to get rid of $t^*,\ s^*$, we have used the
$G$ invariance of the star product.

It is easy to check that the above product agrees with the
multiplication on the crossed product algebra
$C_{c}^{\infty}(M)[[\hbar]]\rtimes G$.
\begin{cor}
The formal deformation quantization of a transformation groupoid
is the crossed product algebra of a formal deformation
quantization of a Poisson manifold  with the $G$ action.
\end{cor}

\begin{rmk}
We may straightforward check that the crossed product algebra is a
formal deformation quantization of the noncommutative Poisson
algebra defined by the crossed Poisson structure in the sense of
Definition \ref{dfn:quant-noncommutative-poisson}. However, our
method in this paper provides a geometrical construction of this
quantization.
\end{rmk}
Furthermore, when $G$ is unimodular, which means $G$ has a
bi-invariant measure and $M$ has a $G$ invariant volume form
$\Omega \index{$\Omega$!volume}$, it is easy to check that
\[
Tr_{\Omega \index{$\Omega$!volume}}(f)=\int d\Omega
\index{$\Omega$!volume}(x) \int d\lambda(\alpha) f(x, \alpha)
\]
defines a trace on $C_c (G, C_c ^{\infty}(M))$ and also on
$C_{c}^{\infty}(M)[[\hbar]]\rtimes G$.
\begin{cor}
If $G$ is unimodular and $M$ has a $G$-invariant measure $\Omega
\index{$\Omega$!volume}$, then $M\times G\rightrightarrows M$ has
a closed deformation quantization with trace $Tr_{\Omega
\index{$\Omega$!volume}}$.
\end{cor}
\begin{rmk}
\label{rmk:morita-eq} If the $G$ action is free and proper, then
$M/G $ is again a Poisson manifold. By the results in
\cite{ri:morita1}, the $C^*-$completion of the crossed product
algebra $C_c^{\infty}(M)\rtimes G$ is Morita equivalent to the
$C^*$-completion of the algebra $C_c^{\infty} (M/G)$ of smooth
compactly supported functions on the quotient. One natural
question is whether this Morita equivalence \index{Morita
equivalence} still holds for the quantized algebras. In the future
work, we will prove that in the case when $G$ is finite and $M$ is
compact, the formal deformation quantizations of $M\rtimes
G\rightrightarrows M$ are Morita equivalent to the corresponding
formal deformation quantizations of $M/G$ in the sense defined by
Bursztyn and Waldman (see \cite{b1:morita} ).
\end{rmk}

\subsection{Strict deformation quantization} In previous sections
of this section, when talking about deformation quantization, we
looked at $\hbar$ as a formal parameter. However, sometimes we do
want $\hbar$ to take some value (even very large). This idea has
inspired Rieffel to introduce a notion of a strict (deformation)
quantization (see Definition \ref{dfn:strict-defor}). In this
section,  we will look at the application of Rieffel's method of
strict deformation quantization \cite{ri:deformation} in strictly
quantizing pseudo \'etale groupoids.\footnote{ Recently, in
\cite{li:strict}, Li has given a ``universal'' way to strictly
quantize Poisson manifolds. His methods can surely be used to
strictly quantize a pseudo Poisson groupoid, but we will not
discuss it here.}

In his lecture notes \cite{ri:deformation}, Rieffel considered a
strict deformation quantization of a Poisson manifold with a
$\reals ^n$ action. Although he was looking at Poisson manifolds,
he was already aware that his method might be able to quantize
some noncommutative Poisson algebras with $\reals ^n$ action.
Notations in the next theorem are explained in the remark which
follows.
\begin{thm}
\label{ri:def}(\cite{ri:deformation}, Theorem 9.3) Let $\alpha$ be
an action of $\reals ^n$ on a $C^*-$algebra $A$, and let $J$ be a
skew-symmetric operator on $\reals ^n$. Let $J$ define a Poisson
bracket, $\{\ ,\ \}$, on $A^{\infty}$. Then $A^{\infty}$ with the
deformed products $\star \index{$\star$} _{\hbar J}$, involutions
$^{*_{\hbar J}}$, and $C^*-$norms $||\ ||_{\hbar}$ as defined in
Chapter 4 of \cite{ri:deformation}, provides a strict deformation
quantization of $A$ in the direction of $\frac {1}{2\pi}\{\ ,\
\}$.
\end{thm}

\begin{rmk}
\begin{enumerate}
\item $A^{\infty}$ stands for the smooth algebra defined by the
$\reals ^n$ action. \item By ``J define a Poisson bracket'', we
mean that for some basis of $\reals ^n$ and the corresponding
derivations $\partial _1,\ \cdots,\ \partial _n$ on $A^{\infty}$
given by the action, the Poisson structure is defined by
\[
\{ f, g\}\stackrel {def}{=}\sum _{j,k}J_{jk}\partial _j
(f)\partial _k (g).
\]
\end{enumerate}
\end{rmk}

In the case of a pseudo \'etale groupoid, if the noncommutative
Poisson structure defined by the formula (\ref{Poisson}) can be
obtained from an $\reals ^n$ action as in Theorem \ref{ri:def}, we
can strictly deformation quantize the groupoid algebra. We
illustrate this by the following example.

\begin{ex}
\label{ex:qdirac} (\cite{tw:dirac})We consider a constant Dirac
structure $\Gamma$ on an $n-$torus $\torus ^n$. As explained in
Example \ref{ex:dirac} of \ref{ex:poisson} , when we choose a
rational transversal subtorus $M$, the reduced foliation
groupoid\footnote{One can construct more than one groupoid from a
foliation. Here, we consider the fundamental groupoid of a
foliation.} is a pseudo Poisson groupoid. Therefore, we can
consider the noncommutative Poisson structure defined on the
groupoid algebra.

The reduced foliation groupoid algebra can be identified with
$C(M)\rtimes \integers ^{k}$. Since $\reals ^{n-k}$ is the
universal covering of $M$, the constant Poisson structure $\pi$ on
$M$ can be lifted onto $\reals ^{n-k}$. From the identification
that $M=\reals ^{n-k}/ \integers ^{n-k}$, $\reals ^{n-k}$ acts on
$C(M)$ by translation commuting with the $\integers ^k$ action;
therefore, there is a well defined $\reals ^{n-k}$ action on the
reduced foliation groupoid algebra. It is not difficult to check
that the noncommutative Poisson structure defined by formula
(\ref{Poisson}) can be obtained from the $\reals ^{n-k}$ action.
Hence,  the conditions of Theorem \ref{ri:def} are satisfied, and
we can strictly deformation quantize the reduced foliation
groupoid algebra. This is an equivalent definition of a
quantization of a constant Dirac structure on an n-torus defined
in \cite{tw:dirac}.
\end{ex}

\end{document}